%% file: revision.tex
\documentclass[onefignum,onetabnum]{siamonline250211} 
\PassOptionsToPackage{hidelinks}{hyperref}
\usepackage{hyperref}

\input{ex_share}

\title{A Joint Variational Framework for Multimodal X-ray Ptychography and Fluorescence Reconstruction}

\begin{document}

\maketitle

\begin{abstract}
Recovering high-resolution structural and compositional information from coherent X-ray measurements involves solving coupled, nonlinear, and ill-posed inverse problems. Ptychography reconstructs a complex transmission function from overlapping diffraction patterns, while X-ray fluorescence provides quantitative, element-specific contrast at lower spatial resolution. We formulate a joint variational framework that integrates these two modalities into a single nonlinear least-squares problem with shared spatial variables. This formulation enforces cross-modal consistency between structural and compositional estimates, improving conditioning and promoting stable convergence. The resulting optimization couples complementary contrast mechanisms (i.e., phase and absorption from ptychography, elemental composition from fluorescence) within a unified inverse model. Numerical experiments on simulated data demonstrate that the joint reconstruction achieves faster convergence, sharper and more quantitative reconstructions, and lower relative error compared with separate inversions. The proposed approach illustrates how multimodal variational formulations can enhance stability, resolution, and interpretability in computational X-ray imaging.  

\end{abstract}



\section{Introduction}

Ptychography has matured into a coherent X-ray microscopy technique in which a focused beam is scanned across the specimen with deliberate overlap while a pixelated detector records far-field diffraction intensities \cite{pfeiffer2018x}. Inverting these intensity-only measurements recovers the specimen’s complex transmission function, enabling quantitative, nanometer-scale imaging across domains from biology to battery materials \cite{zhou2020low, sun2020soft}. Despite its wide applicability, ptychographic reconstruction remains challenging because the inverse problem is nonlinear and nonconvex, with symmetries that can lead to nonuniqueness. Classical solvers, including alternating projections and constraint algorithms \cite{bui2024stochastic,hesse2015proximal,chang2019blind}, ptychographical iterative engine-family (PIE) updates \cite{rodenburg2004phase,maiden2017further}, and maximum-likelihood/variational and Wirtinger-gradient methods \cite{elser2003phase, odstrvcil2018iterative}, can converge rapidly when well initialized, but are susceptible to stagnation and local minima, especially with weak overlap, low counts, or model mismatch.

Complementary contrast is provided by X-ray fluorescence (XRF) \cite{beckhoff2007handbook}, which detects element-specific secondary-emission photons released when incident X-rays excite core electrons. Because each element emits a unique set of characteristic lines, XRF yields quantitative maps of elemental composition and trace distribution, even at very low concentrations. The drawback is that the emitted photons originate from the entire illuminated volume, so the recorded fluorescence map is effectively the convolution of the true elemental distribution with the probe-intensity profile used during scanning. This convolution blurs fine details and constrains the resolution of fluorescence images to the probe size and scan step, even though the data themselves remain chemically specific.

\begin{figure}[H]
    \centering
    \includegraphics[width=0.8\linewidth]{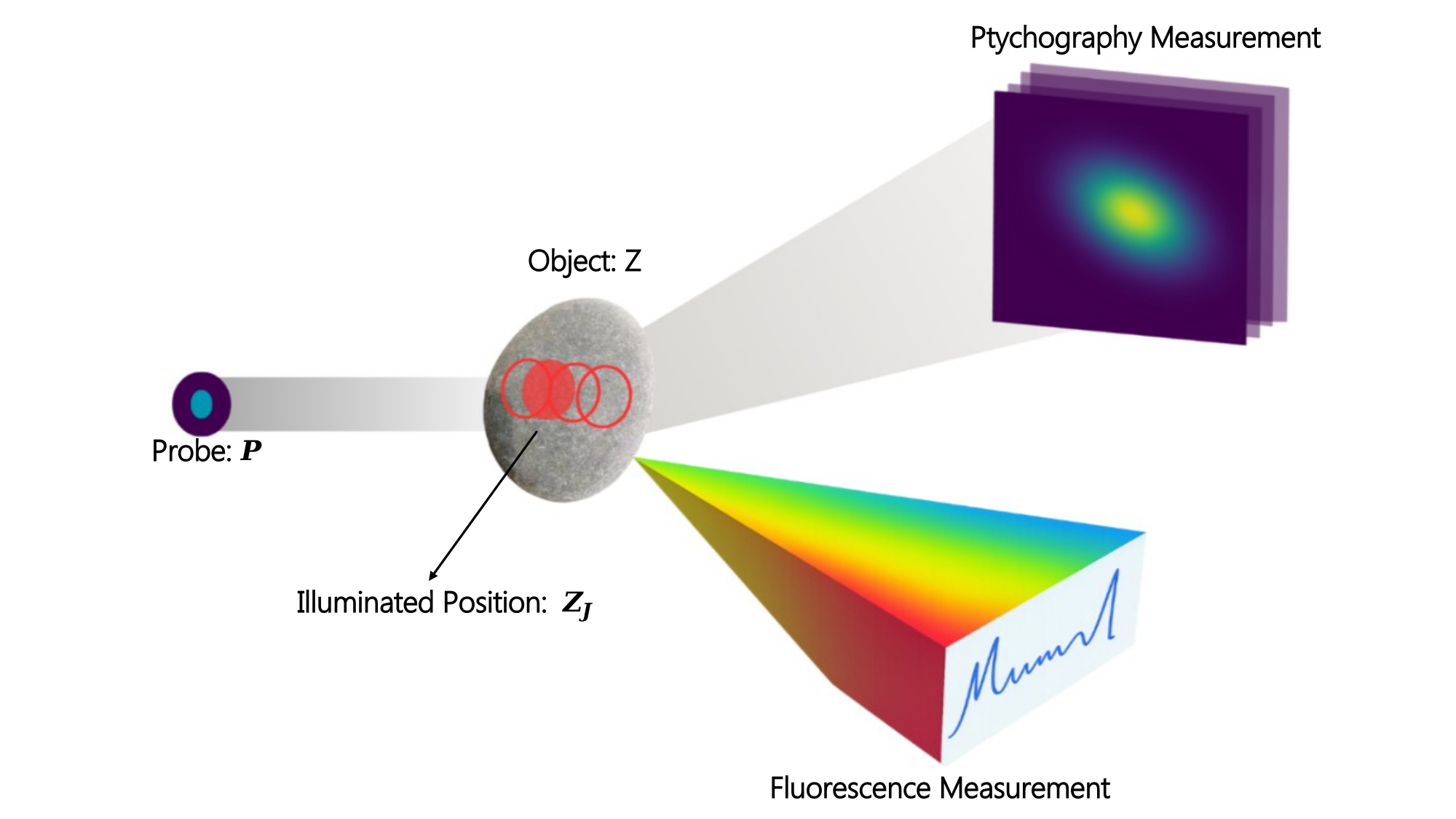}
    \caption{An illustration of a joint ptychographic and fluorescence experiment.}
    \label{fig:joint_experiment}
\end{figure}

{
Ptychography and XRF therefore provide complementary but fundamentally different contrast mechanisms. Ptychography reconstructs the complex transmission function of the specimen: its phase component captures projected refractive-index or electron-density variations, while its amplitude reflects attenuation. This distinction is especially important in weakly absorbing specimens, where absorption contrast may be weak even when phase contrast remains informative. By contrast, XRF does not measure transmission contrast; rather, it records characteristic fluorescence emission and therefore provides chemically specific maps of elemental composition, typically at lower spatial resolution because of probe blur and measurement-response effects. When the two modalities are collected simultaneously (as illustrated in Figure~\ref{fig:joint_experiment}), the ptychographic reconstruction also provides an accurate estimate of the probe, which can then be used as a deconvolution kernel to sharpen the fluorescence signal.}
These complementary strengths naturally motivate a joint inverse formulation that coherently leverages both datasets. 

An early demonstration at synchrotrons combined X-ray ptychography with XRF in a single scan: ptychography supplied a high-resolution complex transmission function while XRF provided element maps; using the reconstructed probe as a point-spread function enabled deconvolution that sharpened fluorescence beyond the optics-limited resolution and contextualized elemental distributions in a biological specimen \cite{vine2012simultaneous,Deng2017}. Di et al.~\cite{di2016optimization} proposed a nonlinear optimization-based joint inversion framework that integrates XRF with X-ray transmission data to reconstruct elemental distributions while simultaneously modeling attenuation and self-absorption, thereby reducing the ill-posedness compared with XRF-only reconstructions. This work established the joint model and solution strategy; subsequent work \cite{Di:17} extended it to joint XRF and transmission tomography, showing improved quantitative accuracy and mitigation of spectral mixing and self-absorption on simulated and experimental datasets. In parallel, the electron-microscopy community advanced closely related ideas. Schwartz et al.~\cite{schwartz2024imaging,schwartz2022imaging} framed fused multimodal electron microscopy as a constrained optimization that links elastic and inelastic scattering, yielding chemical maps with substantially higher signal-to-noise ratio, enabling 10x dose reductions, and providing uncertainty estimates from the Hessian. These works underscore a broader principle of multimodal fusion: tying modalities together through a physics-based objective can trade information across channels to improve quantitative chemical imaging at high spatial resolution.

In this study, we propose a new joint reconstruction framework that explicitly integrates ptychographic and XRF measurements through a {coupling function}. The central idea is to treat the two modalities not as separate post-processing steps, but as coupled observations of the same underlying object: the complex transmission function governs the ptychographic contrast, while the elemental composition drives the fluorescence emission. By embedding these relationships into a single inverse problem, our method leverages the structural sensitivity of ptychography and the elemental specificity of XRF in a mutually reinforcing way. The {coupling function} enforces physical consistency between the modalities, so that each modality regularizes the other within a single optimization framework. The resulting cross-modal coupling improves the conditioning of the otherwise ill-posed inverse problems, suppresses noise and model mismatch, and yields reconstructions with enhanced quantitative accuracy compared with single-modal analyses.

{
The present work uses a lightweight physics-guided coupling between the ptychographic object and the XRF channels to enforce shared spatial consistency, rather than a fully calibrated quantitative XRF transport model. As a result, its accuracy may degrade when important XRF effects, such as calibration errors, detector response, or self-absorption, are not modeled, especially for low-\(Z\) elements. These limitations become particularly important in weak-absorption regimes, where ptychography may remain informative through phase contrast while XRF still detects trace or localized elemental content. In such cases, the formulation should be adapted to fuse high-resolution structural and phase information from ptychography with the chemically specific sensitivity of XRF, rather than assuming a strict one-to-one relation between XRF intensity and transmission absorption.}

{
This paper is therefore best viewed as an initial proof-of-concept study based on simulated data. The numerical experiments are designed to evaluate the behavior of the coupled inverse formulation under controlled conditions, including different probe settings such as in-focus and out-of-focus illumination, rather than to provide full validation on real experimental measurements. Future work will focus on testing the method on experimental data and extending the model to account for calibration uncertainty, detector effects, and other sources of experimental mismatch.}
\section{Joint formulation}
To establish the proposed joint reconstruction framework, we begin by reviewing 
the mathematical models underlying each imaging modality. This provides the 
foundation for formulating a unified loss that couples ptychography and 
fluorescence data.

\subsection{Mathematical background}

In \(2\)D ptychographic phase retrieval, the goal is to reconstruct the complex transmission function of an object \(z = x + y i = |z|e^{i\phi}\in \mathbb{C}^{n\times n} \) from its observed diffraction pattern, where $\phi$ is the phase of the object. 
 The data \(d_j \in \mathbb{R}^{m\times m}\) is collected as {intensity-only diffraction measurements} from the \(\text{j}^{\text{th}}\) scanning position and \(d = \{d_j\}_{j=1}^{N} \in \mathbb{R}^{m\times m\times N}\) is the collection of intensity measurements at all \(N\) scanning positions. A ptychography experiment for one scanning position is modeled by
\[
    d_j = \big| \mathcal{F}(P \odot z_j)\big|^2 + \epsilon_j, 
    \quad j = 1, \cdots, N,
\]
where \(\odot\) denotes {pointwise} matrix multiplication. The operator 
\(\mathcal{F}: \mathbb{C}^{m\times m} \to \mathbb{C}^{m \times m}\) is the 
two-dimensional discrete Fourier transform, and \(|\cdot|^2: \mathbb{C}^{m\times m} 
\to \mathbb{R}^{m\times m}\) computes element-wise intensity. { Given the full object \(z \in \mathbb{C}^{n\times n}\), \(z_j \in \mathbb{C}^{m\times m}\) denotes the local patch illuminated at scan position \(j\). Thus the ptychographic forward model uses the pointwise product \(P \odot z_j\), where \(P \in \mathbb{C}^{m\times m}\) is the probe. } Finally, \(\epsilon_j \in \mathbb{R}^{m\times m}\) models the noise in 
the \(j\)-th measurement.

The overlap ratio quantifies how much the probe’s illuminated footprint at one scan position overlaps the footprint at a neighboring position. Practically, it is reported relative to an effective probe width (e.g., full width at half maximum): a “50\% overlap” means the scan step is about half the probe width along the scan axes. Increasing the overlap increases data redundancy and improves the conditioning of the ptychographic forward operator, which generally enlarges the basin of convergence and makes reconstructions more robust to noise and position/probe errors—albeit at the cost of longer acquisitions and dose. In practice, approximately 60–80\% overlap is commonly used for reliable reconstructions \cite{Maiden:11}.

There are different ways of formulating the ptychographic reconstruction problem \cite{fung2020multigrid}. One such formulation we focus on in this work is the amplitude-based error metric:

\begin{equation} \label{eq:ptychography}
    \min_{P,z} \Phi^{\text{ptyc}}(P,z) =  \frac{1}{2N}\sum_{j = 1}^N \left\| |\mathcal{F}( P \odot z_j) | - \sqrt{d_j}\right\|_{\textit{F}}^2,    
\end{equation}
where the {modulus} operator \(|\cdot|\) makes the optimization problem nonlinear.

On the other hand, we model the expected fluorescence counts for element 
$e$ as a deconvolution with an excitation Point Spread Function (PSF) given by the probe intensity 
\begin{equation} \label{eq:fluorescence}
   \min_{w, P} \; \Phi^{\text{fluor}}(P,w) 
   \;=\; \min_{w,P} \frac{1}{2N_e}  \sum_{e = 1}^{N_e} 
   \left\| \, |P|^2 \ast w_e - D_e \, \right\|_{\textit{F}}^2 ,
\end{equation}
where $\ast$ denotes {two-dimensional discrete convolution on the object grid with kernel $K=|P|^2$, realized by full zero-padded convolution followed by extraction of the central $n\times n$ block. The corresponding adjoint is defined by full convolution with the flipped kernel $\widetilde K(s,t)=K(-s,-t)$, followed by center cropping.} Here, $w_e \in \mathbb{R}^{n \times n}$ is the high-resolution spatial map of element $e$, and $w = [w_e]_{e=1}^{N_e}$ denotes the collection of all $N_e$ XRF elemental maps.  The loss is then computed against 
$D_e \in \mathbb{R}^{n \times n}$, the experimentally measured fluorescence intensity map for element $e$.

{ Throughout the manuscript, we use notation consistently by distinguishing the full object \(z \in \mathbb{C}^{n\times n}\), the local illuminated patch \(z_j \in \mathbb{C}^{m\times m}\), the probe \(P \in \mathbb{C}^{m\times m}\), and the elemental maps \(w_e \in \mathbb{R}^{n\times n}\). Here, \(m\) denotes the probe-window size and \(n\) denotes the full object size, with typically \(m \ll n\).} { We denote the measurements as the collection \(\{d_j\}_{j=1}^N\), where each \(d_j \in \mathbb{R}^{m\times m}\) is the diffraction intensity recorded at scan position \(j\).}

\subsection{Joint reconstruction formulation}
To enhance ptychographic reconstructions under limited data, we integrate fluorescence information by introducing a lightweight, physics-guided coupling between the object's imaginary component \(y\) and its linear attenuation coefficient \(\mu\).
Let \(\mu_e^{E_0}\) be the mass attenuation coefficient at the incident energy \(E_0\).
Under a thin‑sample/single‑scatter model,
\[
\mu \;=\; \sum_{e=1}^{N_e} \mu_e^{E_0}\, w_e.
\]

Physically, attenuation is encoded in the magnitude of the transmission function $|z|=e^{-\frac{1}{2}\int \mu dt}$, with $z=e^{-a}e^{i\phi}$ and $a=-\log|z|$, rather than in $y$ alone. To keep the coupling lightweight, we exploit the global phase gauge of ptychography to apply a constant rotation $z\mapsto \tilde z=e^{-i(\bar\phi+\pi/2)}z$, after which, under a weak-contrast linearization $(|a|\ll 1,\ |\phi-\bar\phi|\ll 1)$, the rotated imaginary part satisfies $\tilde y\approx -e^{-a}\approx -(1-a)$. We therefore treat (a normalized, sign-adjusted) $y$ as a monotone proxy for the line-integrated attenuation, calibrating an affine map from $y$ to $a$ using vacuum/background regions or a thin standard. This approximation avoids log-amplitude transforms and simplifies the joint objective, and is adequate for prototyping algorithmic structure and ablations. Its limitations are clear: the performance degrades when phase variations are large, when the weak-object assumption fails, or for thick/strongly scattering samples.  

Thus, explicitly bridging these two modalities, we introduce a {\emph{coupling function}} that substitutes the absorption term $\mu$ by $y$ in the ptychographic model with the weighted sum of elemental maps:
\[
    y = \sum_{e=1}^{N_e} \mu_e w_e.
\]
This {coupling function} is crucial: it {couples these elemental maps to the imaginary component of the ptychographic object and  ensures that both modalities share a consistent representation of absorption, thereby coupling the elemental concentration maps with the ptychographic reconstruction.} As a result, the two previously independent loss terms, fluorescence data fidelity and ptychographic data fidelity, are unified into a single joint objective. Specifically, we obtain
\[
    \min_{w,P,x}
    \frac{\alpha}{2N_e}\sum_{e = 1}^{N_e} \bigl\|\,|P|^2 \ast w_e - D_e   \bigr\|_{\textit{F}}^2
    + \frac{1}{2N} 
    \sum_{j = 1}^N 
    \bigl\|\,|\mathcal{F}\bigl( P \odot (x + i\sum_{e= 1}^{N_e} \mu_e w_e)_j\bigr) |- \sqrt{d_j}\bigr\|_{\textit{F}}^2,
\]
where $\alpha$ balances the contribution of the two modalities.

For simplicity and identifiability, we restrict this work to the non-blind case, where the probe $P$ is assumed to be known and fixed:

\begin{equation}\label{eq:joint_nonblind}
    \min_{w,x}
    \frac{\alpha}{2N_e}\sum_{e = 1}^{N_e} \bigl\|\,|P|^2 \ast w_e - D_e   \bigr\|_{\textit{F}}^2
    + \frac{1}{2N} 
    \sum_{j = 1}^N 
    \bigl\|\,|\mathcal{F}\bigl( P \odot (x + i\sum_{e= 1}^{N_e} \mu_e w_e)_j\bigr) |- \sqrt{d_j}\bigr\|_{\textit{F}}^2.
\end{equation}

{
The current coupling function neglects several XRF-side effects that are important
for quantitative absorption--composition modeling, including self-absorption,
detector response, spectral overlap and scatter, and geometry-dependent attenuation.
These effects are especially important for low-$Z$ elements, whose low-energy
fluorescence is harder to quantify.  In such regimes, a calibrated attenuation
and transport model would be needed.
}

{
Moreover, when absorption is weak or negligible, the attenuation-based coupling
used in \eqref{eq:joint_nonblind} is no longer the appropriate physical model.
In that case, the joint formulation should instead couple XRF and ptychography
through shared structural information or through a more complete
composition-to-transmission relation, rather than by directly identifying the
elemental maps with the imaginary component of the ptychographic object.
For example, a more complete connection between phase and absorption is provided
by the energy-dependent Kramers--Kronig relation between the dispersive and
absorptive parts of the refractive index, which has been used as an additional
constraint in multi-energy X-ray spectro-ptychography near elemental absorption
edges \cite{hirose2017use}. 
}
{
These extensions are beyond the scope of the present study.  In future work, we
will replace the current proxy with a physically calibrated coupling and, where
needed, a multislice or attenuation-path model.
}
\subsection{Optimal choice of alpha}

We determine \(\alpha\) in \eqref{eq:joint_nonblind} using the GradNorm strategy 
\cite{DBLP:journals/corr/abs-1711-02257}. In GradNorm, the desired gradient norm 
for each task \(i\) is modeled as  
\[
    G_i(z) \;=\; \overline{G(z)} \,\times\, [r_i(z)]^\beta,
\]
where \(G_i\) denotes the gradient norm of task \(i\), \(\overline{G(z)}\) is the 
average gradient norm across all tasks, and \(r_i(z)\) represents the relative 
loss of task \(i\). Under this formulation, the task-specific weights 
\(\alpha_i\) are chosen such that  
\[
    G_i \cdot \alpha_i \;=\; \overline{G(z)},
\]
thereby ensuring that each weighted gradient matches the average gradient norm. Since our joint model involves only two components, we reduce the GradNorm weighting scheme to a single scalar 
parameter \(\alpha\). At initialization, we set  
\[
    \alpha \;=\; \mathcal{O} \left(
\frac{\left\|\,\nabla_{z}\Phi^{\mathrm{ptyc}}(P,z)\big|_{z=z_0}\,\right\|_{2}}
         {\left\|\,\nabla_{w}\Phi^{\mathrm{fluor}}(P,w)\big|_{w=w_0}\,\right\|_{2}} \right),
\]
where \(z_0\) and \(w_0\) are the initial guesses for the objective function~\eqref{eq:joint_nonblind}. This corresponds to GradNorm with \(\beta=1\), enforcing equal 
contributions from the two tasks in terms of gradient norms at the start of training. Balancing the gradient norm is beneficial in the joint optimization because it addresses the loss dynamic mismatch caused by complexity and learning rate imbalance across tasks during optimization. 

In practice, the initial value of \(\alpha\) is highly dependent on the 
experimental configuration, particularly the object size and overlap ratio, 
which affect the relative magnitudes of the ptychography and fluorescence 
gradients. To avoid additional complexity and instability, we fix \(\alpha\) at 
initialization rather than update it dynamically during reconstruction. Also, it is necessary to tune the \(\alpha\) according to the input noise level. For example, a higher weight can be assigned to the fluorescence component in the noiseless setting, while a lower weight is preferable under high-noise conditions.

\section{Spectral and geometric justification of joint reconstruction}

Ptychographic reconstruction is a highly nonlinear and nonconvex inverse problem \cite{chang2021overlappingdomaindecompositionmethods}. Extending to a joint reconstruction framework raises the central challenge of understanding how the coupling between
the ptychography and fluorescence terms reshapes the optimization landscape. To investigate
this, we begin by deriving the Jacobian and gradient of the joint cost function, which characterize the local sensitivity of the residuals to the reconstruction variables. We then analyze
the Hessian, which captures curvature information and governs the local geometry of the loss
surface. Finally, we study the spectrum of the joint Hessian to reveal how incorporating
the fluorescence term modifies the conditioning of the problem, thereby enhancing stability
and convergence. For all analyses, we assume a single-element map \(N_e = 1\) for simplicity;
generalization to multiple-element channels follows directly.

\subsection{Jacobian and gradient derivation}

We first derive the Jacobian expressions and then present the corresponding gradients, which are fundamental to optimization analysis. Here, \((\cdot)^*\) denotes the conjugate transpose, and \(\vec{\cdot}\) denotes vectorization.

Consider the probe matrix 
$P \in \mathbb{C}^{m \times m}$, the object 
$z \in \mathbb{C}^{n \times n}$, and its vectorized form 
$\vec{z} \in \mathbb{C}^{n^2 \times 1}$. We define the residuals
\[
    r_e^{{\text{fluor}}} = |P|^2 \ast w_e - D_e \in \mathbb{R}^{n\times n}, \qquad
    r_j^{\text{ptyc}} = \left | \mathcal{F}(P \odot z_j)\right | - \sqrt{d_j} \in \mathbb{R}^{m\times m},
\]
where $r_e^{{\text{fluor}}}$ is the fluorescence residual and $r_j^{\text{ptyc}}$ is the ptychography 
residual at scan position $j$. Their corresponding contributions to the joint loss \eqref{eq:joint_nonblind} are 
\[
    \Phi_e^{\text{fluor}} = \frac{1}{2 N_e}\|\vec{r_e}^{\text{fluor}}\|_2^2, \qquad 
    \Phi_j^{\text{ptyc}} = \frac{1}{2N}\|\vec{r_j}^{\text{ptyc}}\|_2^2.
\]
For fluorescence, the convolution can be expressed as a linear mapping:
\[
    |P|^2 \ast w_e = \hat{P} \vec{w}_e, \qquad \hat{P} \in \mathbb{R}^{n^2 \times n^2}.
\]

{The vectorized residual vector \( \vec{r_e}^{\text{fluor}}(w_e) \in \mathbb{R}^{n^2 \times 1} \) is differentiated with respect to the vectorized elemental map \( \vec{w}_e \in \mathbb{R}^{n^2 \times 1}\), and we write the corresponding Jacobian and gradient as
\[
    J^{\text{fluor}}_e = \frac{\partial \vec{r_e}^{\text{fluor}}}{\partial \vec{w}_e}=\hat{P} \in \mathbb{R}^{n^2\times n^2}, \qquad 
    \nabla_{\vec{w_e}} \Phi^{\text{fluor}}_e = \frac{1}{N_e}(J^{\text{fluor}}_e)^\top \vec{r_{e}}^{\text{fluor}} \in \mathbb{R}^{n^2 \times 1}.
\]}
In the single-element case \(N_e=1\), this reduces to
\[
    J^{\text{fluor}} = \hat{P}, \qquad
    \nabla_{\vec{w}} \Phi^{\text{fluor}} = (J^{\text{fluor}})^\top \vec{r}^{\text{fluor}}.
\]

For ptychography, let 
$\hat{z}_j = \overrightarrow{\mathcal{F}(P \odot z_j)} \in \mathbb{C}^{m^2 \times 1}$, and let 
$P_j \in \mathbb{C}^{m^2 \times n^2}$ denote the probe operator at position $j$ acting on the vectorized object $z$. 
Let $\mathbf{F} \in \mathbb{C}^{m^2 \times m^2}$ denote the matrix representation of the 2D discrete Fourier transform.

{ 
The vectorized residual $\vec{r_j}^{\rm ptyc}(z)$ is differentiated with respect to the vectorized object variable $\vec{z}$. With the division taken component-wise, we write the corresponding Jacobian as 

\[
    J_j^{\text{ptyc}} = \frac{\partial \vec{r_j}^{\rm ptyc}}{\partial \vec{z}}=\overline{\text{diag}\!\left(\frac{\hat{z}_j}{|\hat{z}_j|}\right)} 
             \,\mathrm{F} P_j \in \mathbb{C}^{m^2 \times n^2}.
\]}
Stacking across all $N$ scan positions gives
\[
    J^{\text{ptyc}} = \begin{bmatrix}
           J_1^{\text{ptyc}} \\
           J_2^{\text{ptyc}} \\
           \vdots \\
           J_{N}^{\text{ptyc}}
         \end{bmatrix}, 
    \qquad
    r^{\text{ptyc}} = \begin{bmatrix}
           r_1^{\text{ptyc}} \\
           r_2^{\text{ptyc}} \\
           \vdots \\
           r_{N}^{\text{ptyc}}
         \end{bmatrix}.
\]
The gradient then follows as
\[
    \nabla_{\vec{z}} \Phi^{\text{ptyc}} = \frac{1}{N} \sum_{j=1}^N (J_j^{\text{ptyc}})^* \vec{r_j}^{\text{ptyc}}
    = \frac{1}{N}\sum_{j = 1}^{N} P_j^* \mathrm{F}^* 
    \text{diag}\!\left(\frac{\hat{z}_j}{|\hat{z}_j|}\right) 
       (|\hat{z}_j| - \sqrt{\vec{d}_j}) \in \mathbb{C}^{n^2 \times 1}.
\]
Finally, { in the single-element case, let \(\theta=(x,w)\) and \(y=\mu w\). Combining both modalities, the Jacobian and gradient of the joint problem \eqref{eq:joint_nonblind} are
\[
J^{\mathrm{joint}}=
\begin{bmatrix}
J^{\mathrm{fluor}}\\
J^{\mathrm{ptyc}}
\end{bmatrix}, \qquad
\nabla_{\theta}\Phi^{\mathrm{joint}}(\theta)=
\begin{bmatrix}
\nabla_{\vec{x}} \Phi^{\mathrm{ptyc}}(x,y)\\
\mu \nabla_{\vec{y}} \Phi^{\mathrm{ptyc}}(x,y)+\alpha \nabla_{\vec{w}}\Phi^{\mathrm{fluor}}(w)
\end{bmatrix}.
\]} These Jacobian expressions highlight how the curvature of the loss surface is directly 
shaped by the structure of $J^{\text{joint}}$, thereby linking spectral analysis to the 
geometric properties of the joint reconstruction problem.

\subsection{Hessian derivation}
The ptychography Hessian takes the block form
\begin{equation} \label{eq:hessian_ptyc}
    H^{\text{ptyc}} = \begin{pmatrix}
        H_{xx} & H_{xw} \\
        H_{wx} & H_{ww}
    \end{pmatrix} \in \mathbb{R}^{2n^2\times 2n^2 }.
\end{equation}
Each block can be expressed in terms of core matrices defined at the \(j^{\text{th}}\) scan position
\[
    C_{1,j} = P_j^* P_j, \quad
    C_{2,j} = \tfrac{1}{m^2} P_j^* \mathrm{F}^* \, \text{diag}\!\left(\tfrac{\sqrt{\vec{d_j}}\odot \hat{z}_j^2}{|\hat{z}_j|^3}\right) \overline{\mathrm{F}} \, \overline{P_j}, \quad
    C_{3,j} = \tfrac{1}{m^2} P_j^* \mathrm{F}^* \, \text{diag}\!\left(\tfrac{\sqrt{\vec{d_j}}}{|\hat{z}_j|}\right) \mathrm{F}P_j,
\]
where \(\overline{\mathrm{F}}\) and \(\overline{P_j}\) denote the conjugate of \(\mathrm{F}\) and \(P_j\), respectively. The four block components of \(H^{\text{ptyc}}\) are then
\[
    H_{xx} = \frac{1}{N} \sum_{j=1}^N \Big[\Re(C_{1,j}) - \tfrac{1}{2}\Re(C_{3,j}) + \tfrac{1}{2}\Re(C_{2,j})\Big],
\]
\[
    H_{ww} = \frac{1}{N} \sum_{j=1}^N \Big[\Re(C_{1,j}) - \tfrac{1}{2}\Re(C_{3,j}) - \tfrac{1}{2}\Re(C_{2,j})\Big],
\]
\[
    H_{xw} = \frac{1}{N} \sum_{j=1}^N \Big[\Im(C_{1,j}) - \tfrac{1}{2}\Im(C_{3,j}) - \tfrac{1}{2}\Im(C_{2,j})\Big], 
    \qquad H_{xw} = H_{wx}^\top.
\]
For the fluorescence component corresponding to one \(w_e\), the Hessian has the form 
\[
    H^{\text{fluor}}_{w_e,w_e} = \nabla_{w_e}^2 \Phi_e^{\text{fluor}} = \frac{1} {N_e}  \frac{\partial}{\partial w_e} (J_e^{\text{fluor}})^\top \vec{r_e}^{\text{fluor}}  = \frac{1}{N_e} \hat{P}^\top \hat{P}.
\]

Since the fluorescence reconstruction is applied to the imaginary part of the object in the ptychography reconstruction, the joint Hessian { is now written with respect to the same stacked variable \(\theta\). In the single-element case, the joint Hessian takes the block form
\begin{equation} \label{eq:hessian_joint}
    H^{\rm joint}(\theta)=
\begin{pmatrix}
H^{\rm ptyc}_{xx} & H^{\rm ptyc}_{xw}\\
H^{\rm ptyc}_{wx} & H^{\rm ptyc}_{ww}+\alpha H^{\rm fluor}_{ww}
\end{pmatrix},
\end{equation}
so that the fluorescence term contributes only to the \(w\)-\(w\) block, while the cross-terms arise from the ptychographic dependence on \(y=\mu w\).}

This block structure makes explicit that the fluorescence term adds positive curvature to the ptychography Hessian. To quantify how this changes optimization, we analyze the eigenvalue distribution of the joint Hessian under a deliberately challenging synthetic setup: object size { $ n = 62$, one element map ($N_e=1$), probe size $ m = 16$,} $6.25\%$ overlap, noiseless data, and a Hessian of size $7688\times 7688$. The very low overlap makes the baseline (ptychography-only) problem highly ill-conditioned. Because the Hessian spectrum encodes local curvature of the loss, it offers direct insight into optimization behavior \cite{liao2021hessianeigenspectrarealisticnonlinear}. For shape comparison, we normalize eigenvalues by the maximum, $\tilde\lambda_i=\lambda_i/\lambda_{\max}(H)$.

\begin{figure}
  \centering
  \includegraphics[width=0.8\linewidth]{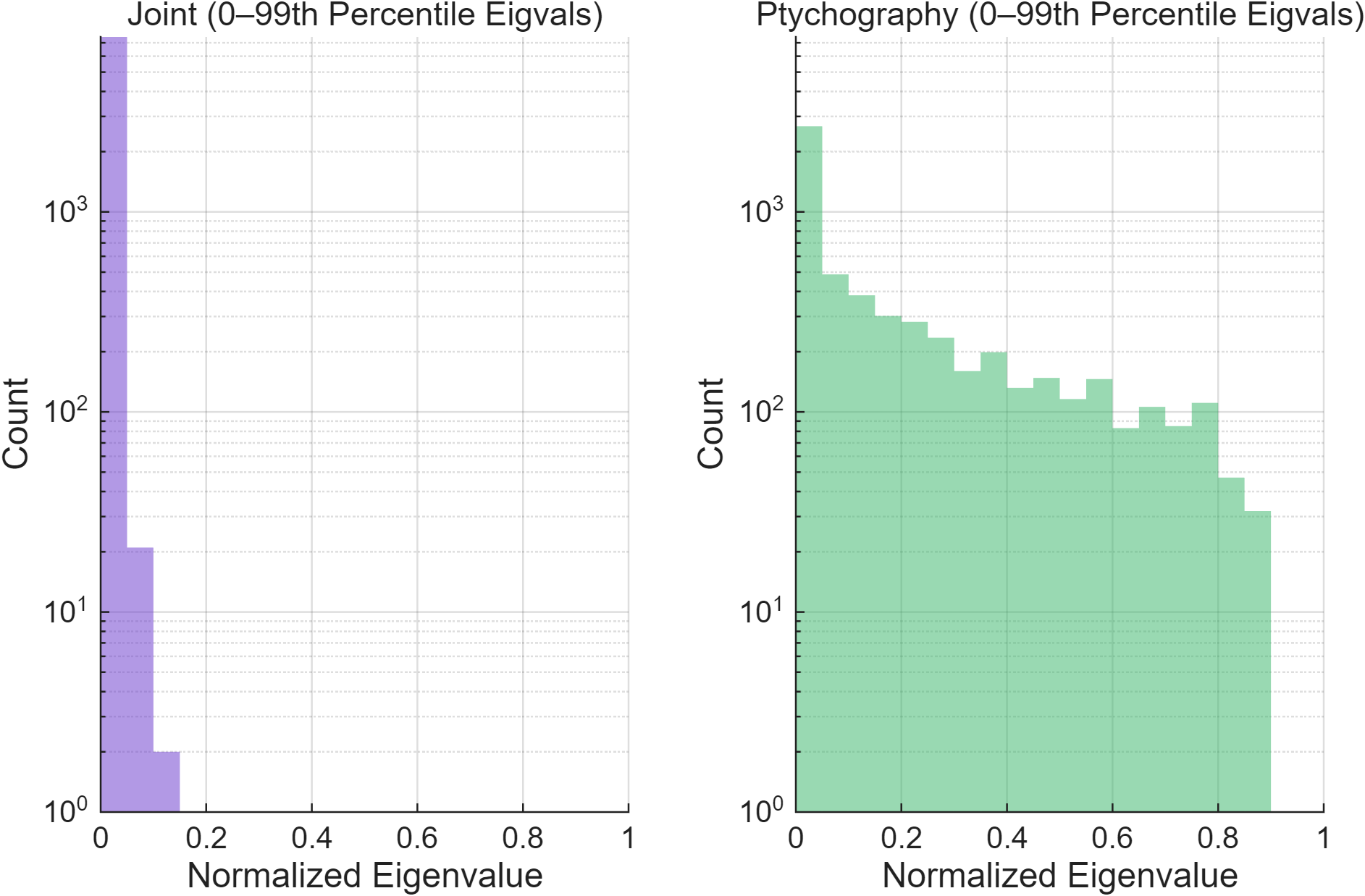}
  \caption{Normalized eigenvalue distribution comparison in reconstruction.}
  \label{fig:eig_dist_comparison}
\end{figure}

Figure~\ref{fig:eig_dist_comparison} shows that the ptychography-only Hessian has a relatively even distribution of eigenvalues across $[0,1]$, reflecting many weakly curved directions. By contrast, the joint Hessian exhibits a two-phase spectrum: the bulk ($\approx 99\%$) is tightly clustered near zero, while a small fraction ($\approx 1\%$) form outliers at much larger values, generating a clear spectral gap. This structure indicates that the joint objective has stronger curvature along a few informative directions, enhancing local identifiability and improving the effectiveness of descent updates within these subspaces \cite{sagun2017eigenvalueshessiandeeplearning}.

\begin{figure}
  \centering
  \includegraphics[width=0.8\linewidth]{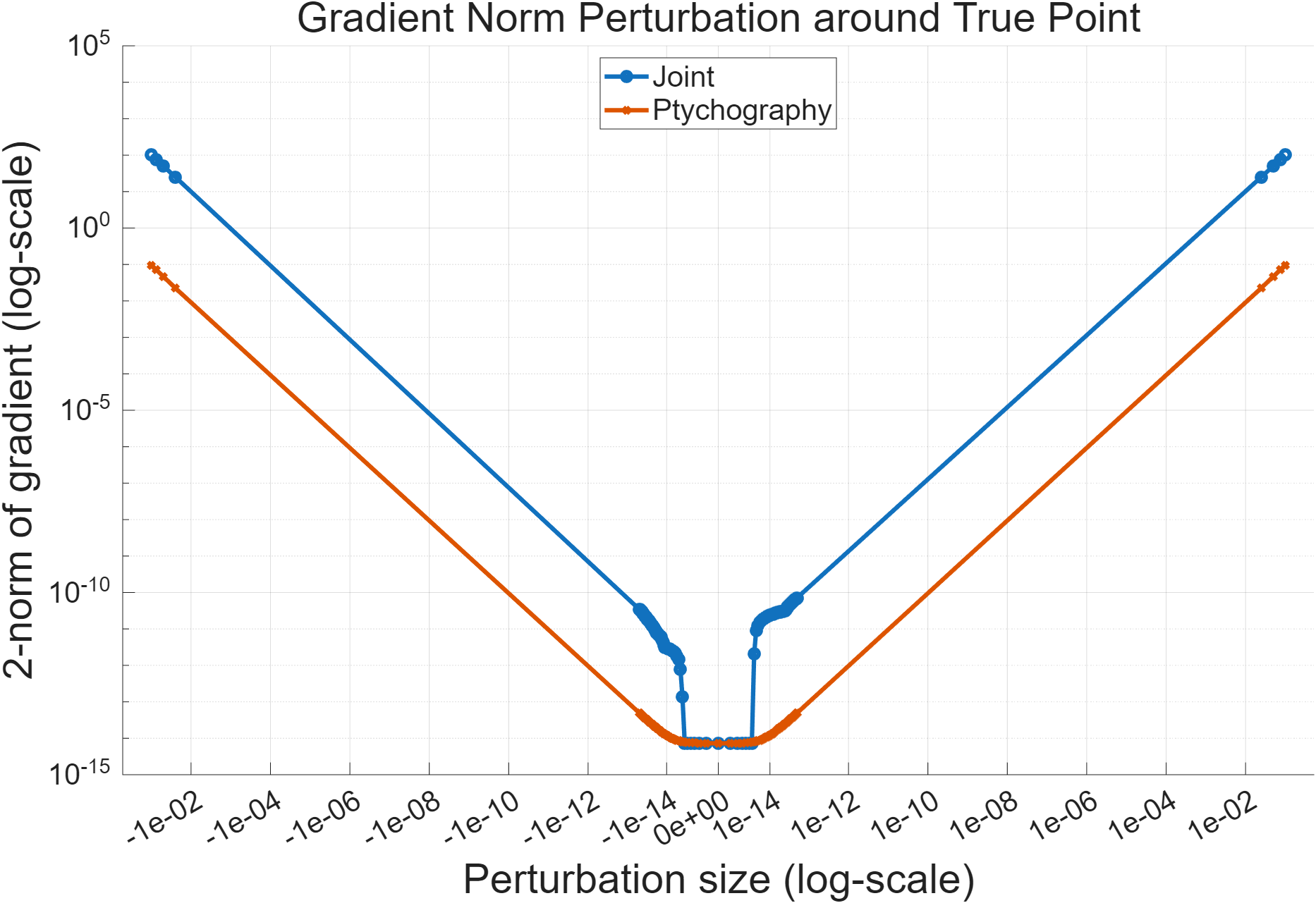}
  \caption{True point gradient norm perturbation log-scale plot with illustration.}
  \label{fig:gradient_norm_true_point}
\end{figure}

A complementary gradient–perturbation experiment around the ground truth (Figure~\ref{fig:gradient_norm_true_point}) reinforces this interpretation. For perturbations $\Delta$ taken along the eigenvector corresponding to the largest eigenvalue, we vary the perturbation magnitude $|\Delta|$ continuously up to $0.1$ in both positive and negative directions. In this neighborhood, the local linearization $\nabla\Phi(z^\star+\Delta)\approx H\Delta$ holds for a loss function $\Phi$ evaluated near the global optimum $z^\star$. The faster growth of $|\nabla\Phi^{\text{joint}}(z^\star+\Delta)|$ relative to $|\nabla\Phi^{\text{ptyc}}(z^\star+\Delta)|$ within the $10^{-14}$ perturbation range reflects the larger eigenvalues in the informative subspace of $H^{\text{joint}}$. In this regime, stronger curvature tightens the link between loss reduction and reconstruction error.

\subsection{Evolving structure of the Hessian}

Having analyzed the joint Hessian’s spectrum at a fixed point, we now ask whether the joint formulation activates more \emph{informative curvature directions} throughout optimization. While the spectrum remains dominated by a small set of outliers, the joint model also lifts additional modes above a fixed tolerance. To quantify this effect, we track the numerical rank
\[
r_\tau(H) := \#\{\sigma_i(H)\ge \tau\}, \qquad \tau=10^{-4},
\]
for both the ptychography Hessian in Eq.~\eqref{eq:hessian_ptyc} and the joint Hessian in Eq.~\eqref{eq:hessian_joint} at several iterations.

Table~\ref{tab:joint_tolerance_rank} reports the results. Across all iterations, $r_\tau(H^{\text{joint}})$ consistently exceeds $r_\tau(H^{\text{ptyc}})$. The relative gain is modest ($\approx 2.6$–$4.3\%$) in early iterations but grows to $10.1$–$10.7\%$ at later stages.

In high-dimensional models, numerical rank serves as a proxy for the dimensionality of informative curvature \cite{ambartsumyan2020hierarchicalmatrixapproximationshessians}. A low rank indicates redundancy, with many parameters contributing little to curvature, whereas a higher rank means more directions are independently active. From this perspective, the joint Hessian not only produces strong outliers (anchoring the landscape with sharp curvature in a few highly informative directions) but also broadens the subspace above the noise floor by engaging additional modes. Thus, the fluorescence term contributes in two complementary ways: it sharpens critical features through large eigenvalues, and it expands the set of active directions through rank gains. Together, these effects reduce redundancy and make optimization more stable and effective.

\begin{table}[H]
  \centering
  \begin{tabular}{lcccc}
    \toprule
    \textbf{Iteration} & $r_\tau(H^{\text{ptyc}})$ & $r_\tau(H^{\text{joint}})$ & $\Delta$ & \%$\uparrow$ \\
    \midrule
    $25$  & $7201$ & $7390$ & $+189$ & $+2.6\%$ \\
    $50$  & $6867$ & $7163$ & $+296$ & $+4.3\%$ \\
    $150$ & $6045$ & $6656$ & $+611$ & $+10.1\%$ \\
    $250$ & $5982$ & $6619$ & $+637$ & $+10.6\%$ \\
    $375$ & $5972$ & $6622$ & $+650$ & $+10.9\%$ \\
    $500$ & \textbf{5978} & \textbf{6621} & \textbf{+643} & \textbf{+10.7\%} \\
    \bottomrule
  \end{tabular}
  \caption{Numerical rank $r_\tau(H)$ with $\tau=10^{-4}$ (number of singular values $\ge \tau$). 
  The joint Hessian consistently achieves a higher rank, indicating more active optimization directions.}
\label{tab:joint_tolerance_rank}
\end{table}

\subsection{Sharpness and smoothness}

Building on the outlier–bulk spectrum and rank growth observed in the previous subsections, 
we next probe the geometry of the loss surface through two-dimensional slices. 
As shown in \cite{li2018visualizinglosslandscapeneural}, highly irregular landscapes can impede optimization.

We visualize the loss around the current reconstruction $\theta^\star$ by perturbing along two orthonormal directions $v_1,v_2$:
\[
    f(\alpha,\beta) \;=\; \Phi\!\bigl(\theta^\star + \alpha v_1 + \beta v_2\bigr),
\]
with $v_1,v_2$ normalized following the filter-wise strategy of \cite{li2018visualizinglosslandscapeneural}. 
To contrast sharp versus flat subspaces, we choose $v_1,v_2$ from Hessian eigenvectors: 
for sharp directions, the eigenvectors corresponding to the $200$th and $201$st largest eigenvalues; 
for flat directions, those associated with the $4000$th and $4001$st eigenvalues. 
For comparability, we fix the normalized perturbation range $(\alpha,\beta)\in[-5,5]^2$ and plot normalized values 
$f(\alpha,\beta)/f_{\max}$ with $f_{\max}=\max_{(\alpha,\beta)} f(\alpha,\beta)$.

\begin{figure}[t]
    \centering
    \includegraphics[width=0.8\linewidth]{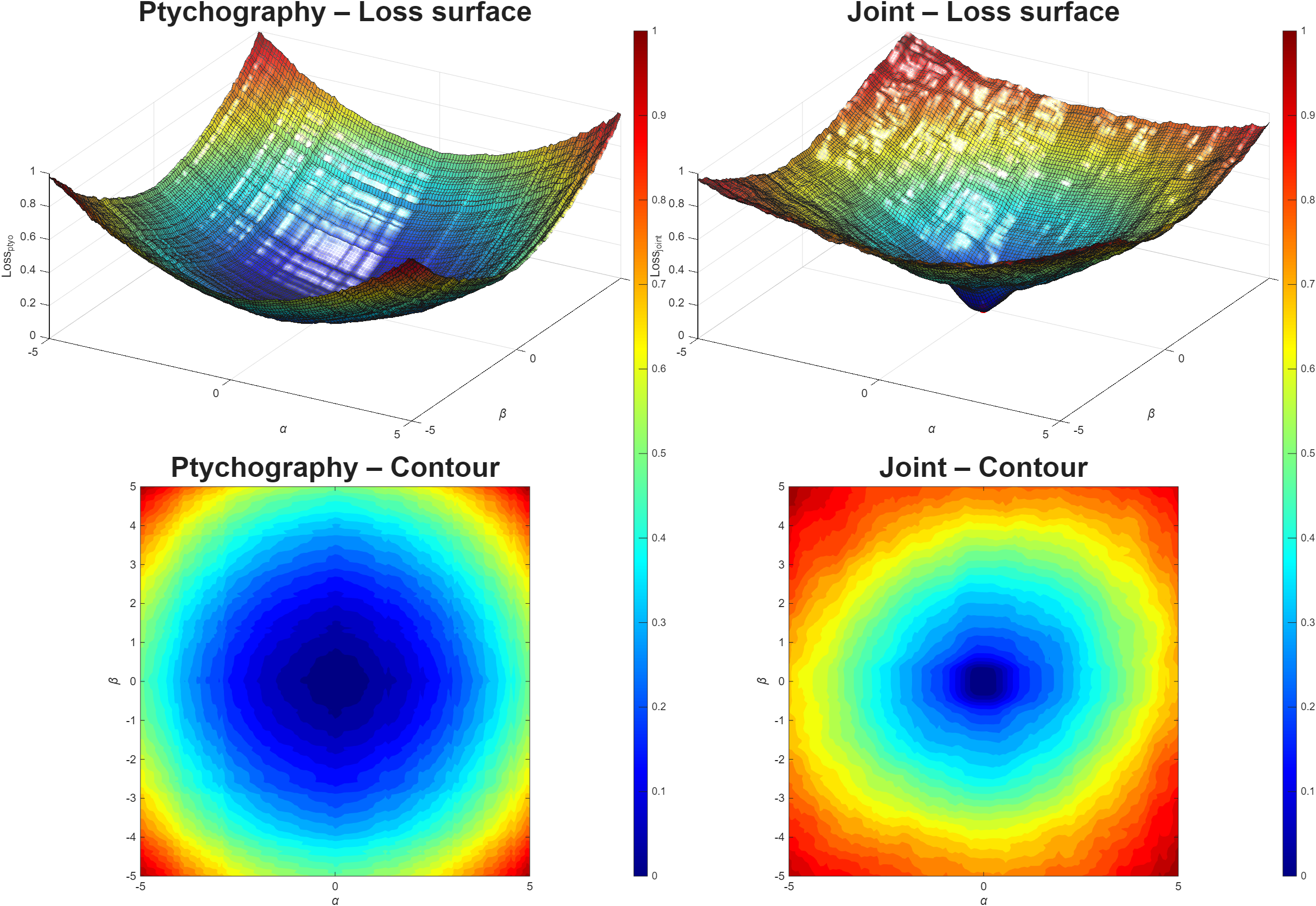}
    \caption{Loss surface around $\theta^\star$ along two sharp directions (eigenvectors for the $200$th and $201$st largest eigenvalues). Comparison between ptychography (left) and joint (right). The joint landscape exhibits a steeper, more decisive basin.}
    \label{fig:loss_surface_eig200}
\end{figure}

\begin{figure}[t]
    \centering
    \includegraphics[width=0.8\linewidth]{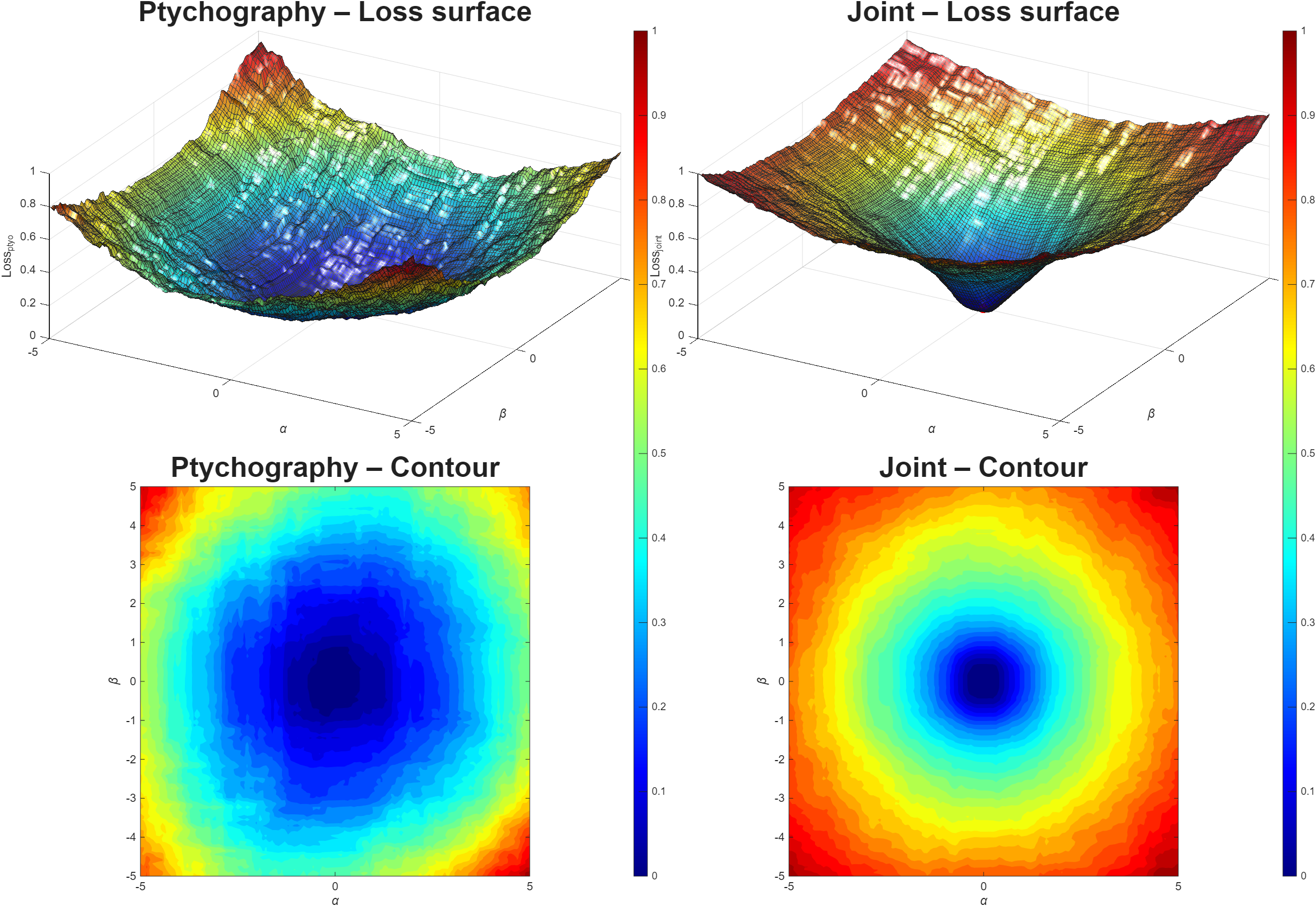}
    \caption{Loss surface around $\theta^\star$ along two flat directions (eigenvectors for the $4000$th and $4001$st largest eigenvalues). Comparison between ptychography (left) and joint (right). The joint landscape is smoother and more convex-like.}
    \label{fig:loss_surface_eig6000}
\end{figure}

The results (Figs.~\ref{fig:loss_surface_eig200}, \ref{fig:loss_surface_eig6000}) show that the joint formulation 
produces a markedly sharper local basin along the leading-eigenvector plane, while also yielding a smoother, more convex-like landscape along flatter directions. 
This combination is advantageous: strong curvature in a small informative subspace anchors critical features, 
while a broader set of directions with non-negligible curvature reduces redundancy and enables progress in regions that are otherwise nearly flat. 
Together, these effects shape a loss surface that is both globally easier to traverse  
and locally more decisive, resulting in improved stability and convergence.

\section{Numerical experiment}

In this section, we systematically compare ptychographic, XRF and joint reconstruction methods under varying experimental conditions to assess the practical benefits of our proposed formulation. All reconstructions are performed using the truncated Newton method \cite{nash2000survey}, a large-scale optimization algorithm that efficiently handles nonlinear problems through inexact Hessian approximations and robust line-search updates. 

We evaluate performance across five settings: varying overlap ratio, noise robustness, scalability with object size, handling of multiple elemental maps, and robustness with different sample features. All code is implemented in MATLAB \footnote{Code repository: \url{https://github.com/EricZou007/Joint-Reconstruction.git}} and executed on an AMD \(\text{R}9-8945\) CPU. 

\begin{figure}
    \centering
    \includegraphics[width=0.8\linewidth]{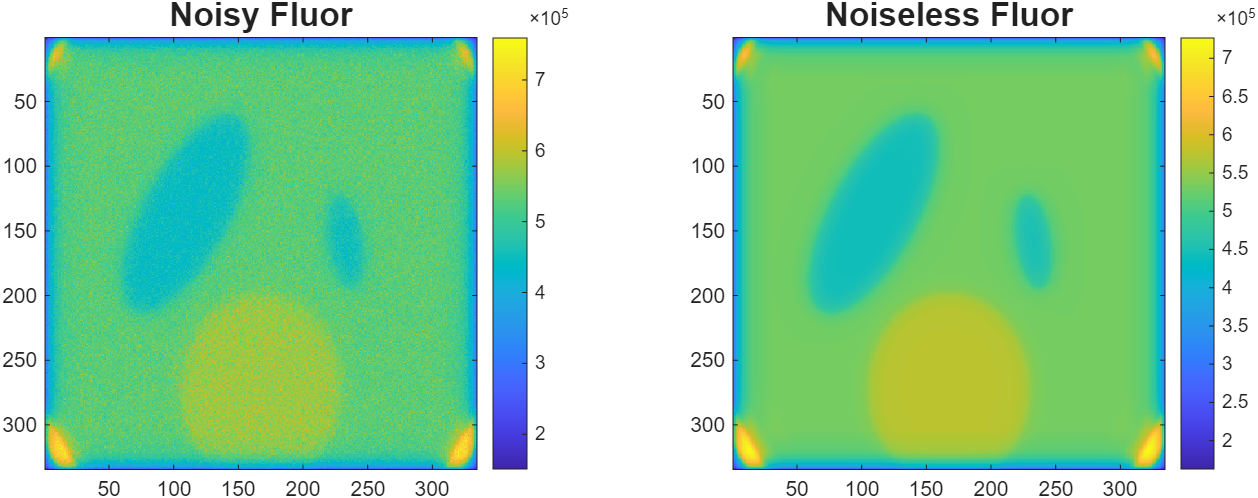}
    \caption{{Noisy XRF map (left) and noise-free XRF map (right).}}
    \label{fig:fluor_measurement}
\end{figure}

\paragraph{Dataset}
We generate synthetic ptychography and fluorescence measurements from known complex objects and elemental concentrations. {The simulated fluorescence measurements are generated by convolving the high-resolution elemental map with the probe intensity profile, so the degree of XRF blurring is determined by the probe-intensity profile. To make this explicit, we  show representative examples of the blurred XRF elemental map in Figure \ref{fig:fluor_measurement}.}

\begin{figure}
  \centering
  \begin{subfigure}[t]{0.5\linewidth}
    \centering
    \includegraphics[width=0.7\linewidth]{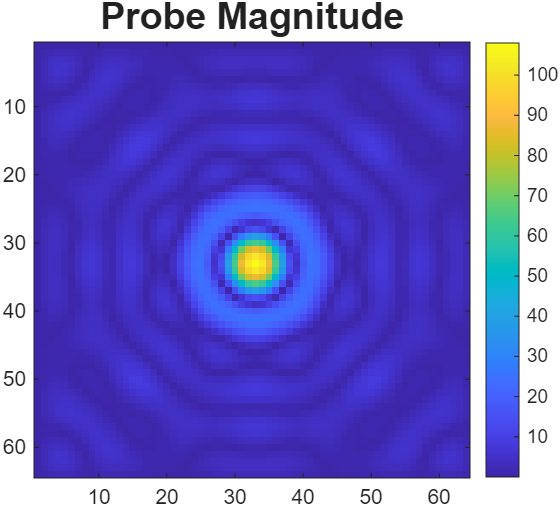}
  \end{subfigure}%
  \hfill
  \begin{subfigure}[t]{0.5\linewidth}
    \centering
    \includegraphics[width=0.7\linewidth]{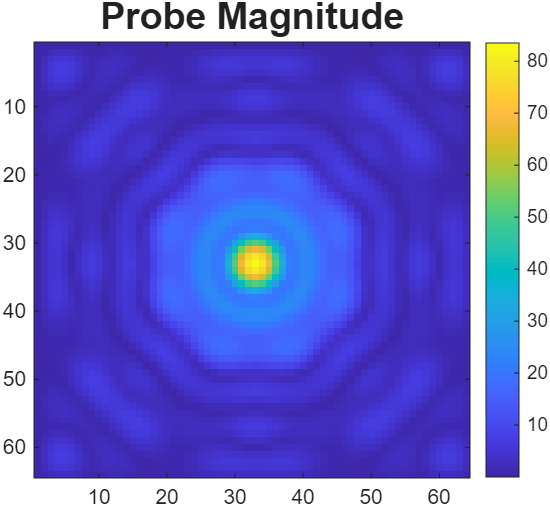}
  \end{subfigure}
  \caption{{In-focus (left) and out-of-focus (right) simulated probe magnitude used in the numerical experiments.} }
  \label{fig:probe_focus}
\end{figure}

\paragraph{Probe and scanning setup}
The probe is simulated as a Fresnel zone plate following \cite{Rodenburg2008PtychographyAR}. Scanning positions form a regular grid, corresponding to a \(50\%\) overlap ratio by default. Increased overlap improves identifiability and reduces ill-posedness in ptychography.  { Figure~\ref{fig:probe_focus} shows the probe on the same pixel grid as the reconstructed object, with $m=64$.}  It is worth noting the probe’s periphery away from the center has near-zero amplitude, which can induce ill-conditioning.

{
Different probe conditions (i.e., in-focus versus out-of-focus) can significantly affect the conditioning of the reconstruction problem. In the in-focus setting, the probe intensity is more concentrated near the center, which reduces measurement diversity. In the default setting, our numerical experiments use an out-of-focus probe. }

\paragraph{Noise model}
To simulate photon-counting noise, we inject Poisson noise into both ptychographic and fluorescence measurements. { Poisson noise is used as a first-order model because both modalities record photon counts at the detector, and in the photon-limited regime the dominant uncertainty is shot noise. This is a standard assumption in photon-counting X-ray imaging and is consistent with common XRF and ptychographic noise models.} 

Given a noiseless intensity \(\hat{d}_j\), we define the noisy data and quantify the noise level as
\[
d_{j,k} = \eta \cdot \text{Poisson}\left(\frac{\hat{d}_{j,k}}{\eta}\right), \qquad \text{Noise Level} = \sqrt{\frac{\eta}{\operatorname{mean}(\hat{d}_j)}} \times 100 \%
\]
where \(\eta\) is a dimensionless parameter controlling the noise level. A larger \(\eta\) implies fewer photons and higher relative noise. We set \(\text{Noise Level} = 3\%\) by default.

\paragraph{Evaluation metrics}
We report two complementary metrics:
\begin{itemize}
    \item \textbf{Loss:} The objective value \(\Phi(z)\) quantifies the consistency between the reconstruction and the measured data. For a fair comparison, we evaluate the joint model by computing the ptychographic loss function \(\Phi^{\text{ptyc}}\) and fluorescence loss function \(\Phi^{\text{fluor}}\) using the solution \(z^{\text{joint}}\) obtained from the joint formulation~\eqref{eq:joint_nonblind}.  

    \item \textbf{Reconstruction error:} The mean squared error (MSE) between the reconstructed complex object $z$ and the ground truth $z^*$: 
    \[
        \frac{1}{n^2} \left\|z - z^*\right\|_F^2.
    \]
\end{itemize}
This distinction is essential because in ill-posed problems, a small residual may not imply a correct solution, due to the phenomenon of semi-convergence in the presence of noise \cite{10.5555/1805888}. {
Here, the fluorescence-only reconstruction refers to the stand-alone solution of the XRF subproblem in Equation~\ref{eq:fluorescence}, i.e., a non-blind deconvolution with a known excitation PSF given by the probe intensity. Under the present simplified coupling, this reconstruction estimates only the attenuation-linked component used in the joint formulation (i.e., the imaginary part). Accordingly, the XRF-only baseline is compared only against this attenuation-linked quantity, rather than against the full complex object, since fluorescence data alone do not determine the ptychographic phase. For consistency, we therefore evaluate the mean squared error (MSE) of the imaginary component when comparing fluorescence, ptychography, and joint reconstruction results.
}


\paragraph{Stopping criterion}
To ensure a fair comparison, all methods and cases run for \(100\) optimization steps. Final results report both loss and MSE at termination.

\subsection{Effect of overlap ratio on performance}

We first investigate how the overlap ratio affects reconstruction quality. 
We consider two representative cases: 83\% overlap and 20\% overlap. 

\begin{figure}
  \centering
  \includegraphics[width=0.8\linewidth]{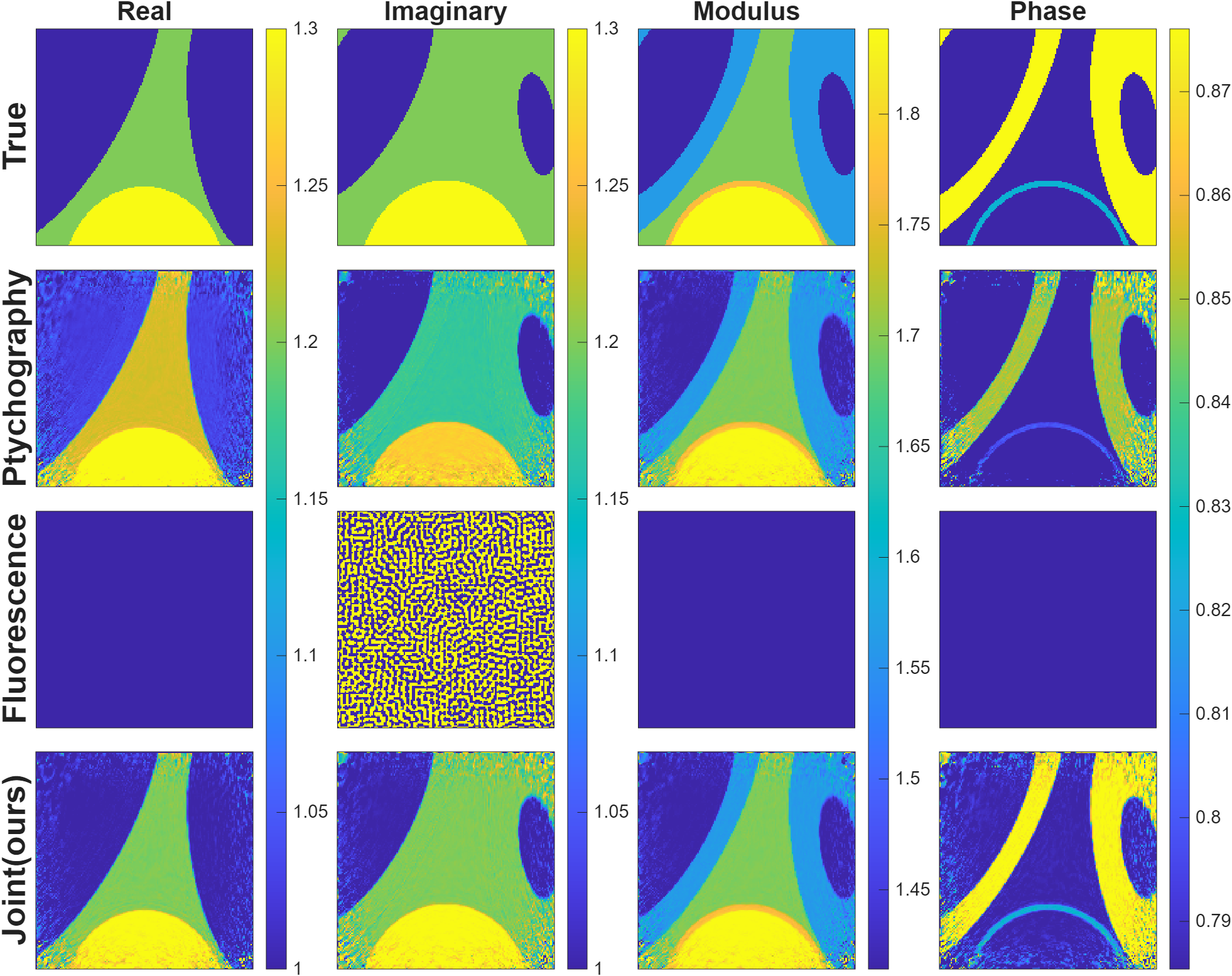}
  \caption{Reconstruction results for ptychography only, fluorescence only and the joint method on a synthetic object of size $n=154$ with probe size $m=64$, using $\text{Noise Level} = 3\%$ and overlap ratio $0.83$. The joint method achieves more accurate reconstructions.}
  \label{fig:recon_obj154}
\end{figure}

\begin{figure}
  \centering
  \includegraphics[width=0.8\linewidth]{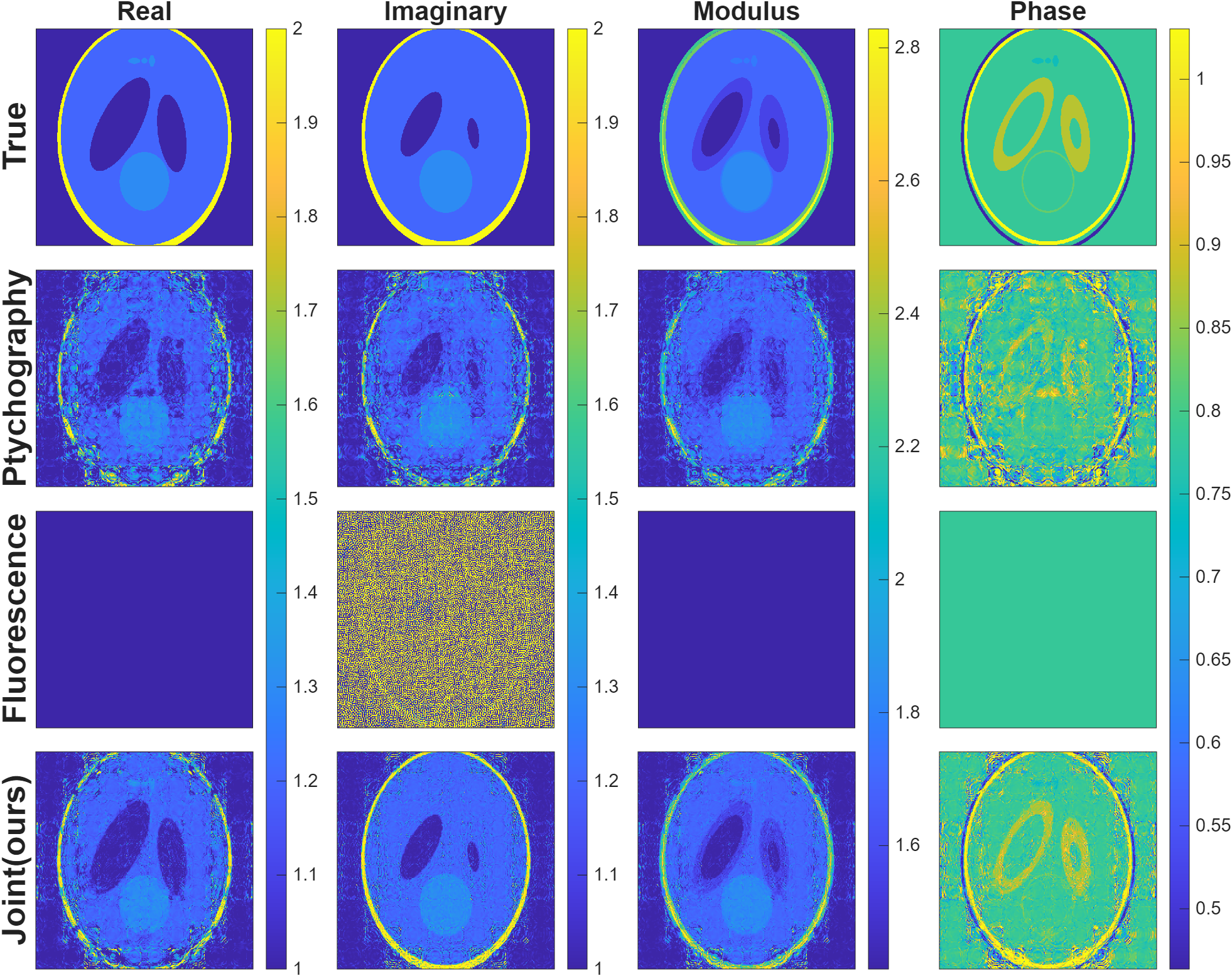}
  \caption{Reconstruction results for ptychography and the joint method on a larger synthetic object of size $n=514$ with probe size $m=64$, with \(\text{Noise Level} = 3\%\) and 20\% overlap. The joint method outperforms ptychography, especially in recovering boundary details of the imaginary component.}
  \label{fig:recon_obj514}
\end{figure}

For the high-overlap case (83\%), we observe a significant gap between the MSE of the ptychography reconstruction and the joint reconstruction. For the low-overlap case (20\%), again, we see an 
{enlarged} gap between the MSE for ptychography and the joint. {This improvement is caused by the additional regularization provided by the fluorescence modality.} Representative reconstructions for the true object, ptychography, fluorescence, and the joint method are shown in Figures~\ref{fig:recon_obj154} and \ref{fig:recon_obj514}.  In both cases, especially in the low-overlap case, the joint approach (augmented with fluorescence data) achieves higher-quality reconstruction of both components more consistent with the true image. 

{By contrast, in the absence of regularization, the {fluorescence-only} baseline reflects a challenging deconvolution problem involving a broad excitation probe and limited measurement diversity, rendering the inversion highly sensitive to noise and numerical instability.} This instability is further amplified by the square-root operation. In both cases, the fluorescence method exhibited semi-convergence to noise, yielding a blurred, noisy image after a limited number of iterations. {The benefit of the joint formulation is that the ptychographic term suppresses deconvolution artifacts by enforcing cross-modal consistency.}

\begin{figure}
  \centering
  \begin{subfigure}[t]{0.45\linewidth}
    \centering
    \includegraphics[width=0.8\linewidth]{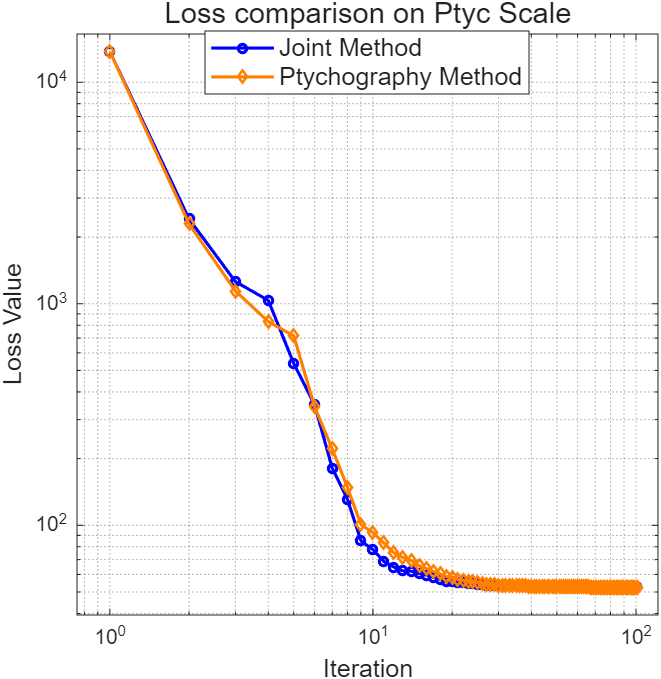}
  \end{subfigure}%
  \hfill
  \begin{subfigure}[t]{0.45\linewidth}
    \centering
    \includegraphics[width=0.8\linewidth]{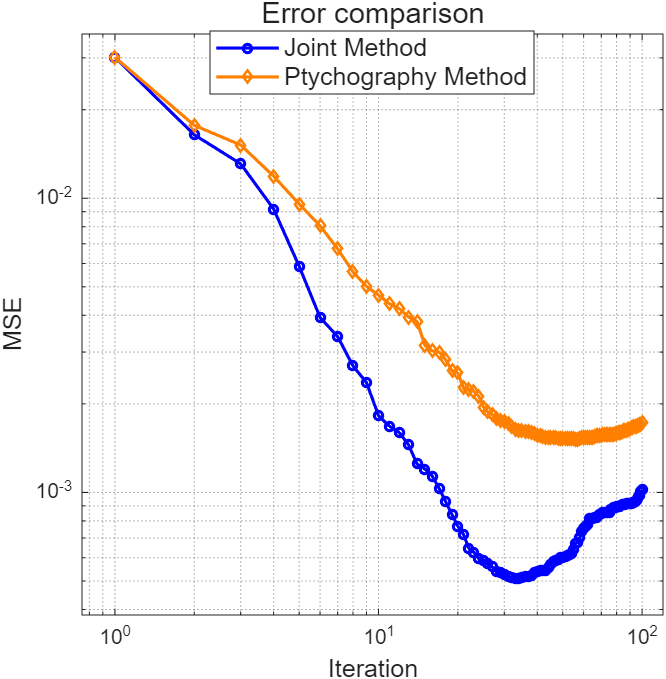}
  \end{subfigure}
  \caption{Log–log plots of loss (left) and reconstruction error (right) for ptychography and the joint method on the synthetic object shown in Figure \ref{fig:recon_obj154} with probe size $m=64$, with \(\text{Noise Level} = 3\%\) and 83\% overlap. The joint method achieves comparable convergence while attaining lower reconstruction error.}
  \label{fig:loss_error_overlap0.75}
\end{figure}

\begin{figure}
  \centering
  \begin{subfigure}[t]{0.45\linewidth}
    \centering
    \includegraphics[width=0.8\linewidth]{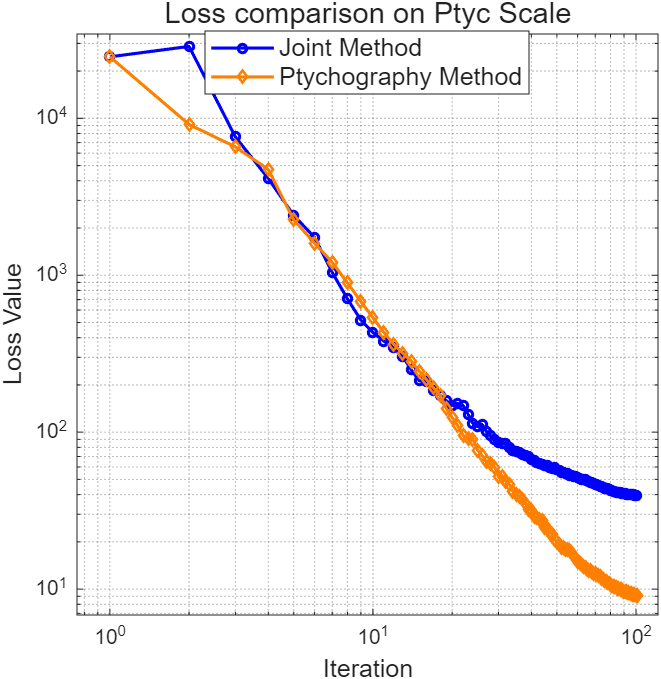}
  \end{subfigure}%
  \hfill
  \begin{subfigure}[t]{0.45\linewidth}
    \centering
    \includegraphics[width=0.8\linewidth]{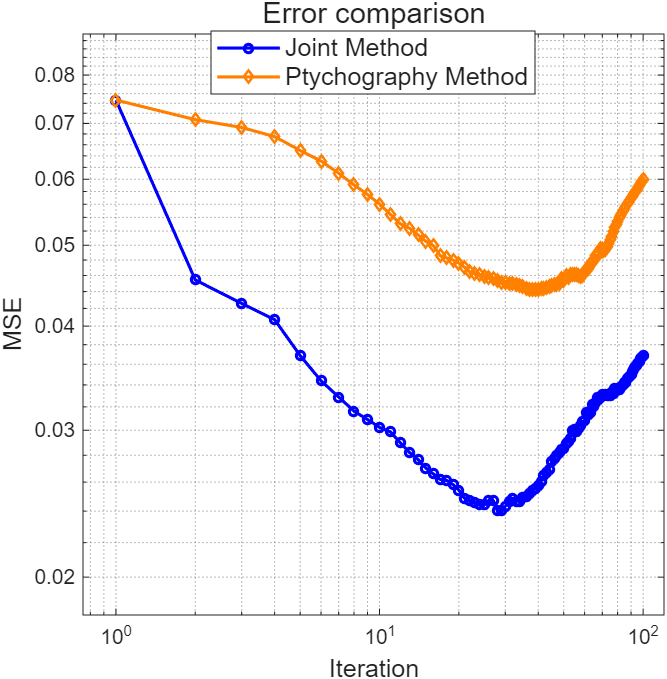}
  \end{subfigure}
  \caption{Log–log plots of loss (left) and reconstruction error (right) for ptychography and the joint method on the synthetic object shown in Figure \ref{fig:recon_obj514} with probe size $m=64$, with $\text{Noise Level} = 3\%$ and 20\% overlap. The joint method maintains lower error throughout optimization compared with ptychography-only method.}
  \label{fig:loss_error_overlap0.2}
\end{figure}

Figures~\ref{fig:loss_error_overlap0.75} and \ref{fig:loss_error_overlap0.2} report the corresponding log–log plots of loss and reconstruction error. The trends confirm our earlier analysis: the joint objective introduces additional curvature in the imaginary component, broadening the effective search directions during optimization. As a result, the joint method attains substantially lower MSE at the same loss level, highlighting its robustness across overlap conditions.

\begin{figure}
  \centering
  \begin{subfigure}[t]{0.45\linewidth}
    \centering
    \includegraphics[width=0.8\linewidth]{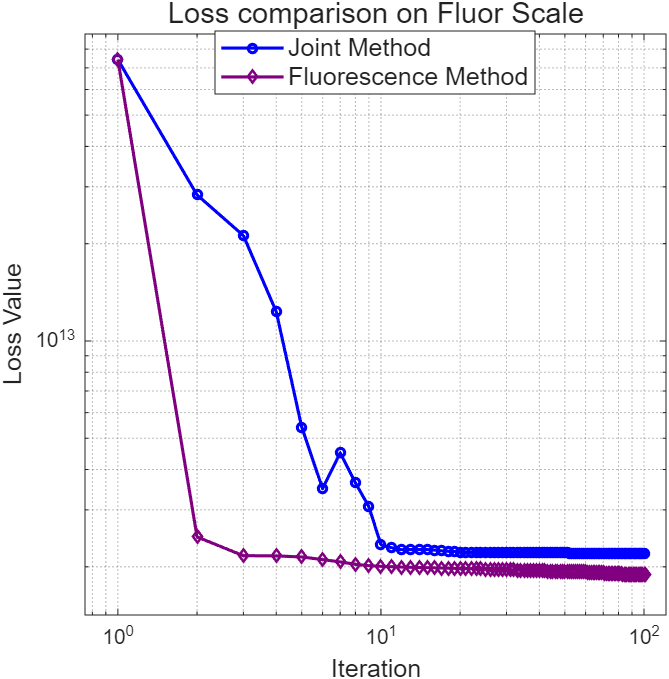}
    \vspace{0.6em}
    \includegraphics[width=0.8\linewidth]{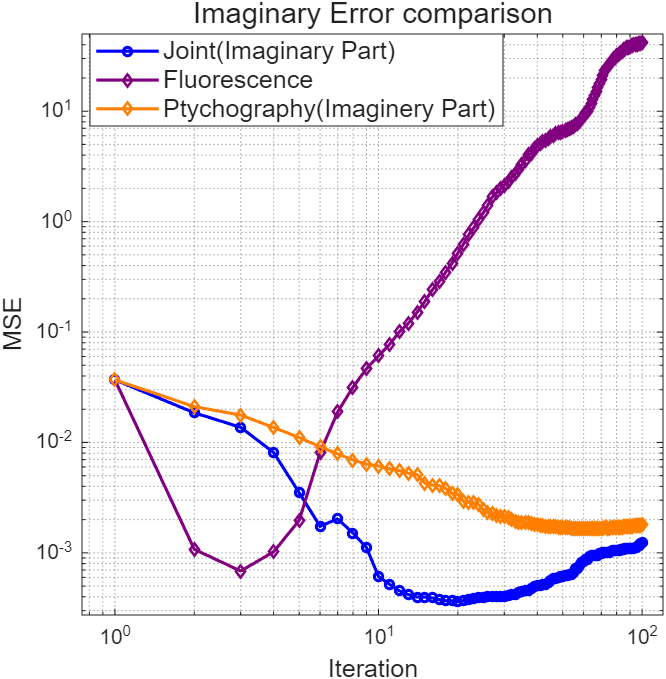}
    \caption*{$n=154$, overlap $=0.83$}
  \end{subfigure}%
  \hfill
  \begin{subfigure}[t]{0.45\linewidth}
    \centering
    \includegraphics[width=0.8\linewidth]{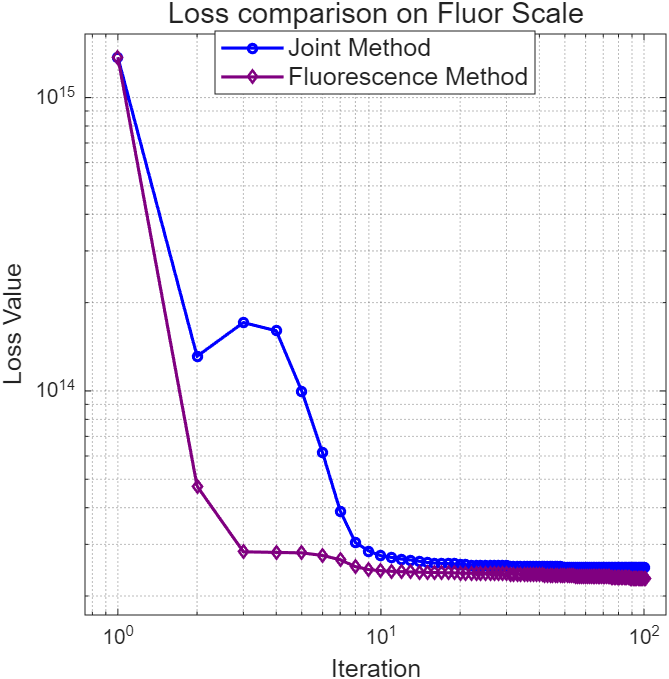}
    \vspace{0.6em}
    \includegraphics[width=0.8\linewidth]{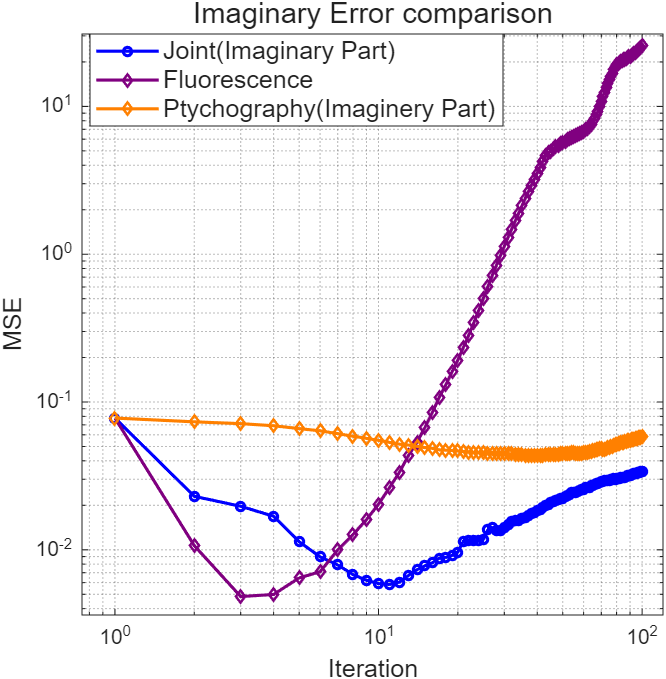}
    \caption*{$n=514$, overlap $=0.2$}
  \end{subfigure}

  \caption{Log--log plots of objective loss (top) and reconstruction error (bottom) for fluorescence and the joint method under $\text{Noise Level} = 3\%$ noise and \(m = 64\) using the synthetic objects shown in Figures \ref{fig:recon_obj154} (left) and \ref{fig:recon_obj514} (right). Joint method exhibits stronger robustness to the input noise while maintaining lower MSE.  }
  \label{fig:loss_fluor_error_overlap_ratio}
\end{figure}

As shown in Figure \ref{fig:loss_fluor_error_overlap_ratio}, the fluorescence MSE only decreases for the first few iterations and then rises due to the semi-convergence effect. By incorporating constraints from the ptychography data, the joint reconstruction demonstrates improved robustness against Poisson noise while maintaining relatively lower MSE, consistent with our theoretical expectation.

\subsection{Stability under noise}

The joint method maintains high accuracy even in the presence of noise. To evaluate stability, we compare the joint and ptychography methods across different noise levels by adding Poisson noise to both ptychography and fluorescence data. Two cases are considered: \(\text{Noise Level} = \{0\%, 10\%\}\).  

\begin{figure}
  \centering
  \begin{subfigure}[t]{0.45\linewidth}
    \centering
    \includegraphics[width=\linewidth]{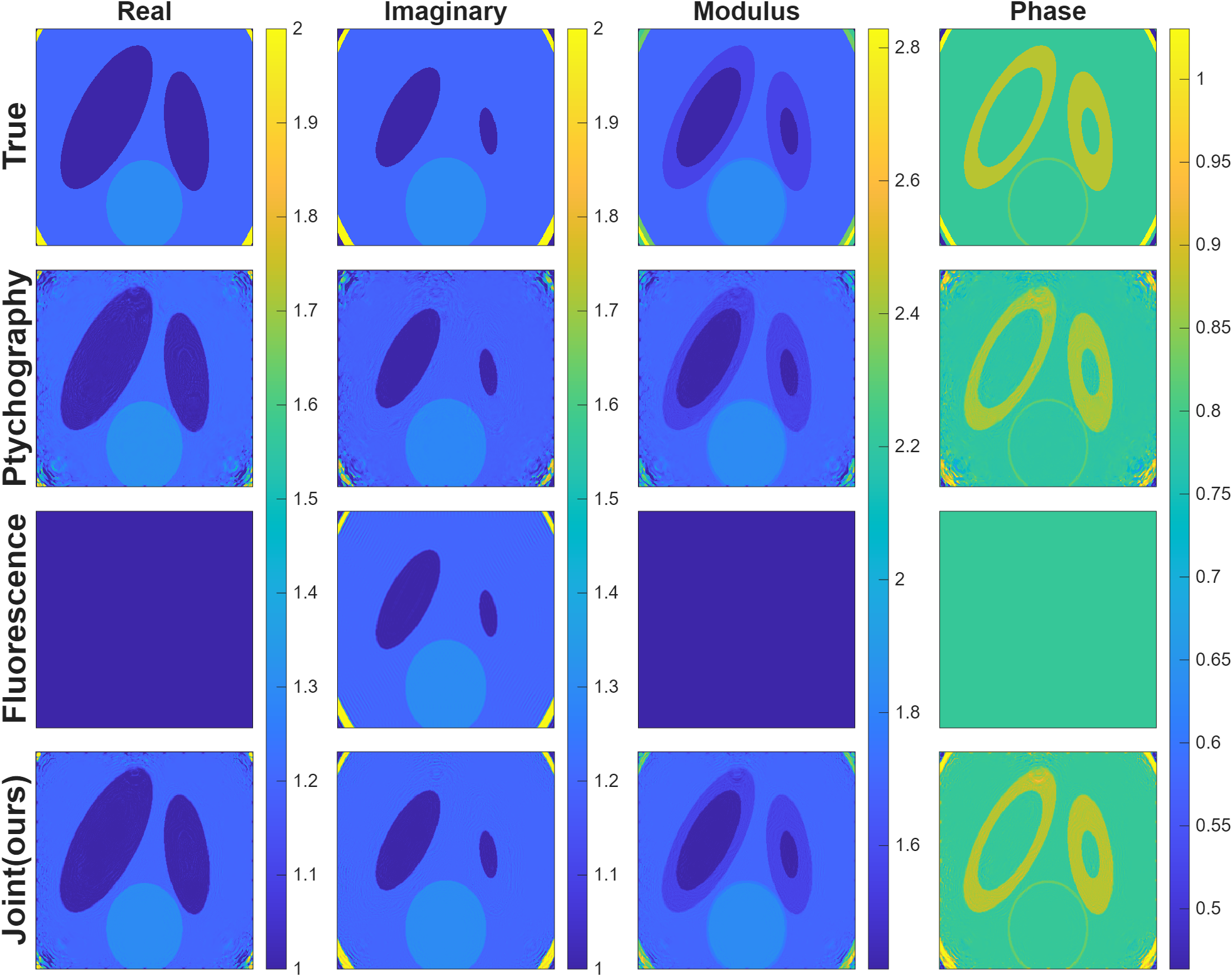}
    \caption*{\(\text{Noise Level} = 0\%\)}
  \end{subfigure}%
  \hfill
  \begin{subfigure}[t]{0.45\linewidth}
    \centering
    \includegraphics[width=\linewidth]{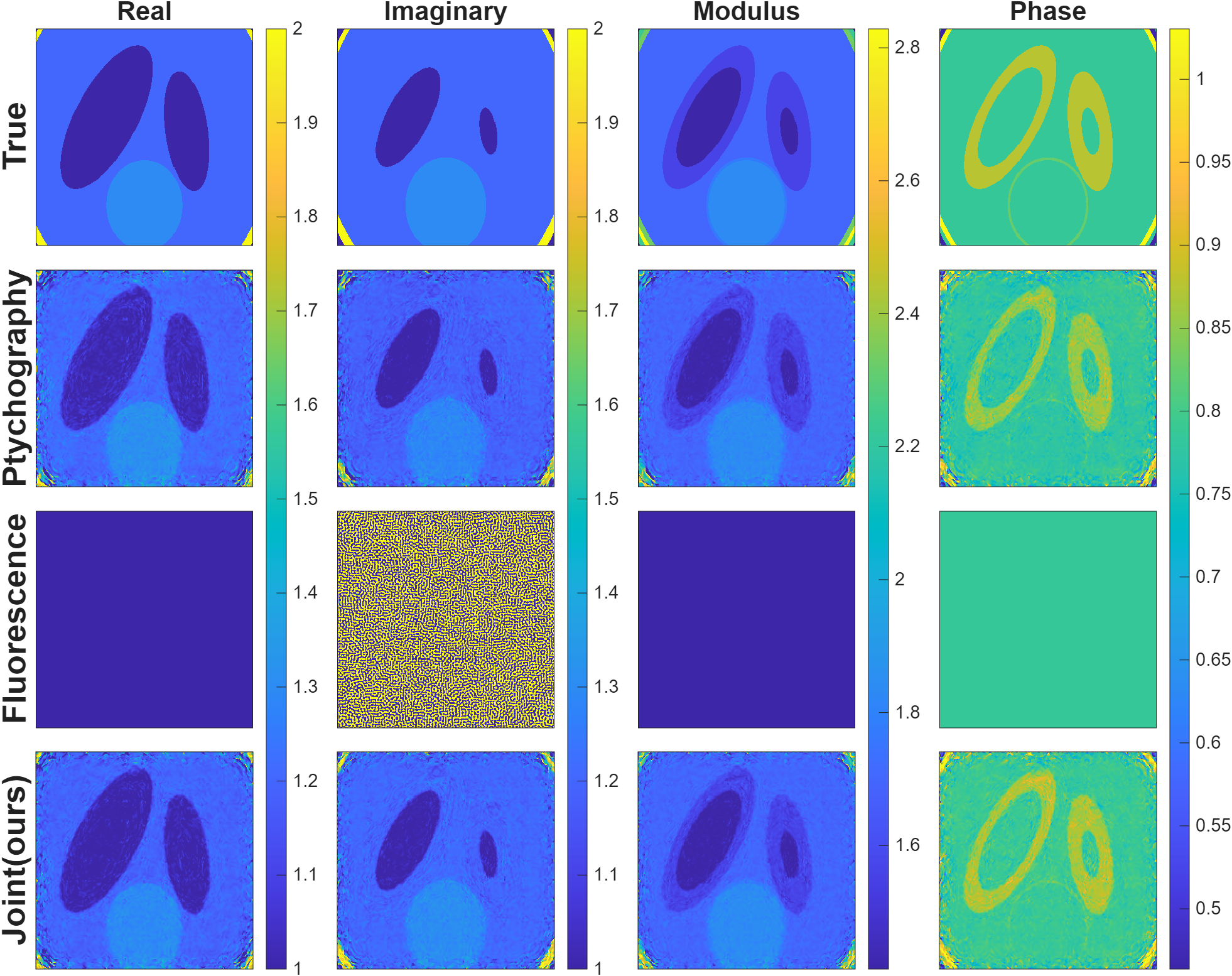}
    \caption*{\(\text{Noise Level} = 10\%\)}
  \end{subfigure}

  \caption{Reconstruction results for ptychography and joint method applied to a synthetic object with \(n=334\), 53\% overlap, and probe size \(m=64\), under two noise settings.}
  \label{fig:noise_recon}
\end{figure}

\begin{figure}
  \centering
  \begin{subfigure}[t]{0.45\linewidth}
    \centering
    \includegraphics[width=0.8\linewidth]{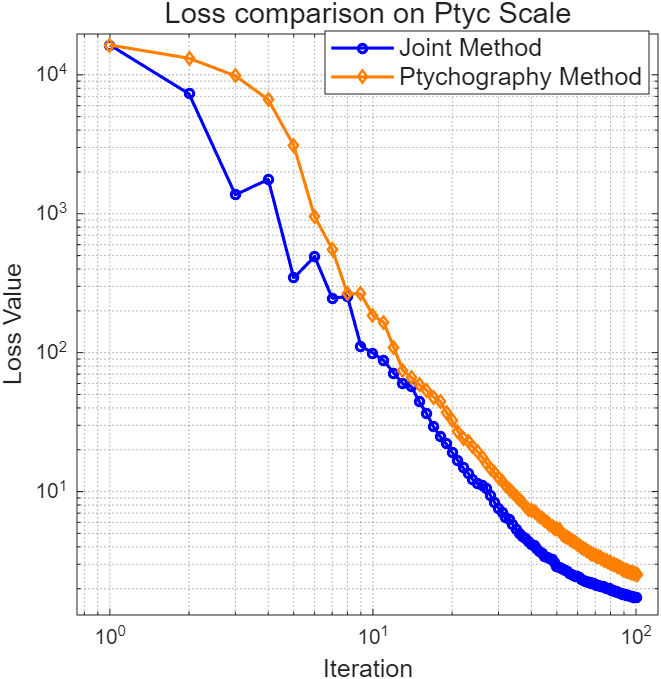}
    \vspace{0.6em}
    \includegraphics[width=0.8\linewidth]{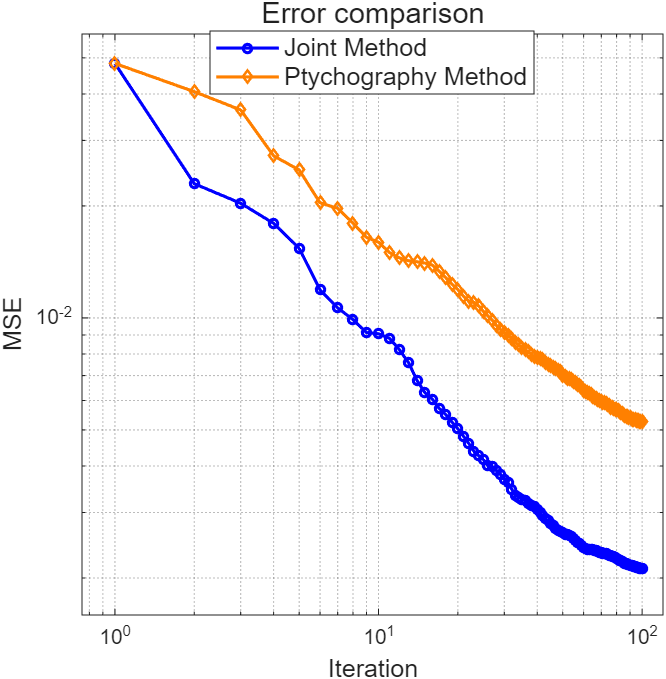}
    \caption*{\(\text{Noise Level} = 0\%\): loss (top), error (bottom)}
  \end{subfigure}%
  \hfill
  \begin{subfigure}[t]{0.45\linewidth}
    \centering
    \includegraphics[width=0.8\linewidth]{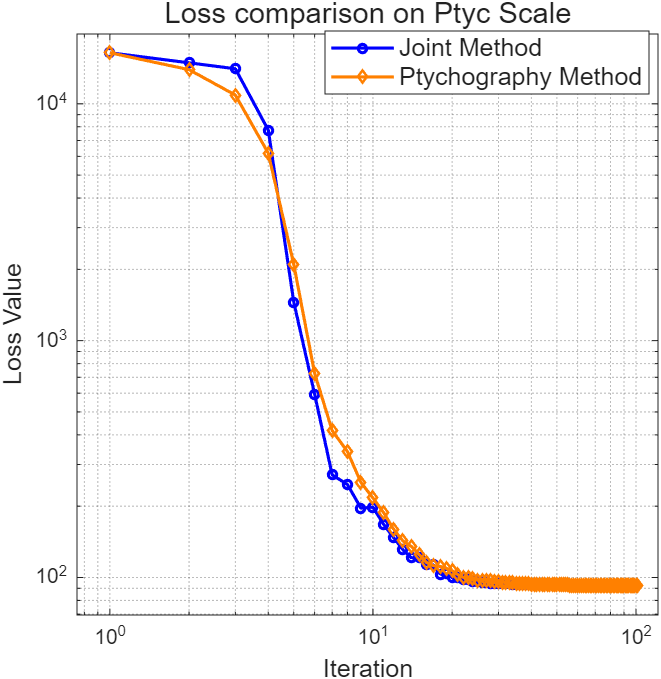}
    \vspace{0.6em}
    \includegraphics[width=0.8\linewidth]{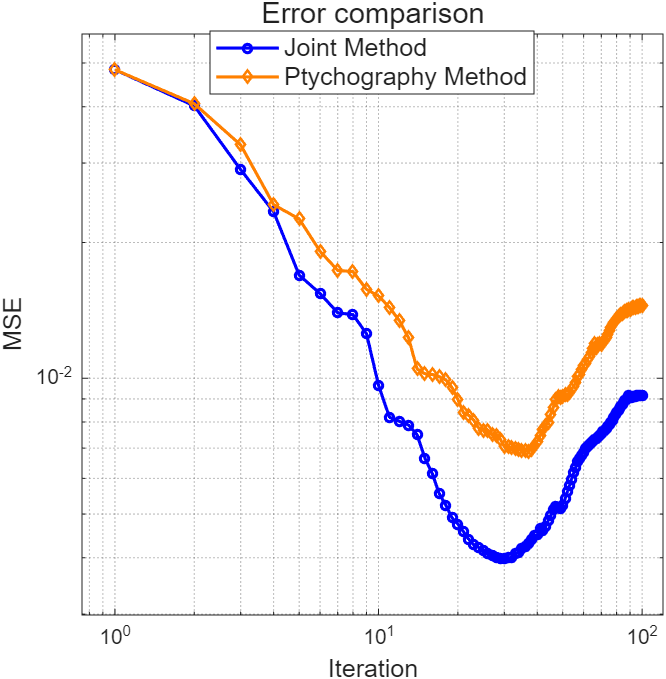}
    \caption*{\(\text{Noise Level} = 10\%\): loss (top), error (bottom)}
  \end{subfigure}

  \caption{Log--log plots of loss (top) and reconstruction error (bottom) for ptychography and joint method applied to the synthetic object shown in Figure \ref{fig:noise_recon} with \(n=334\), probe size \(m=64\), and 53\% overlap.}
  \label{fig:loss_errro_noise_compare}
\end{figure}

\begin{figure}
  \centering
  \begin{subfigure}[t]{0.45\linewidth}
    \centering
    \includegraphics[width=0.8\linewidth]{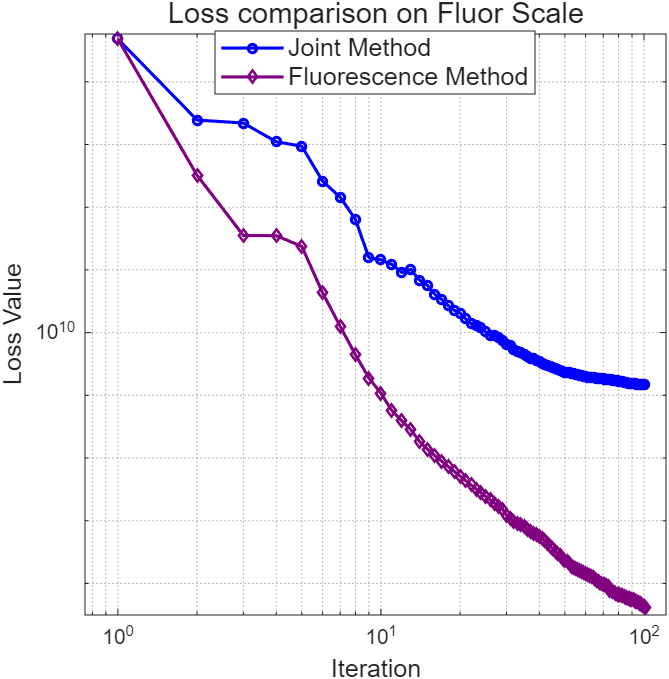}
    \vspace{0.6em}
    \includegraphics[width=0.8\linewidth]{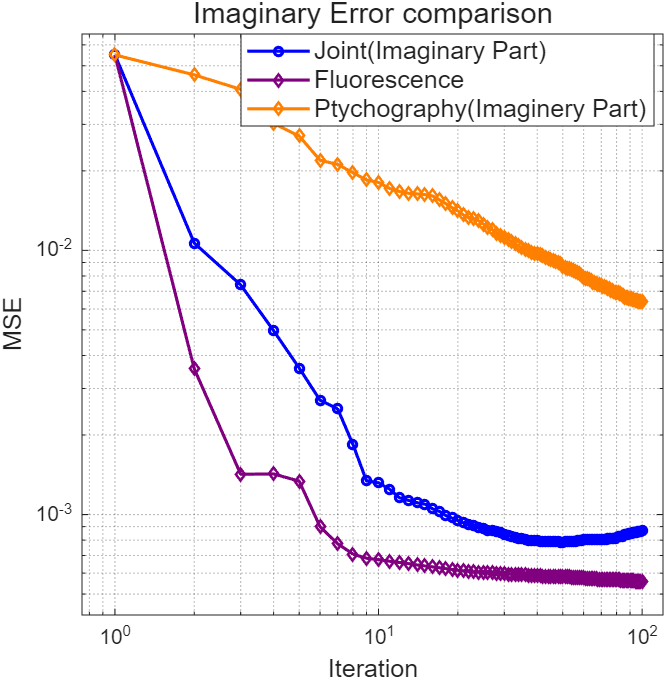}
    \caption*{\(\text{Noise Level} = 0\%\)}
  \end{subfigure}%
  \hfill
  \begin{subfigure}[t]{0.45\linewidth}
    \centering
    \includegraphics[width=0.8\linewidth]{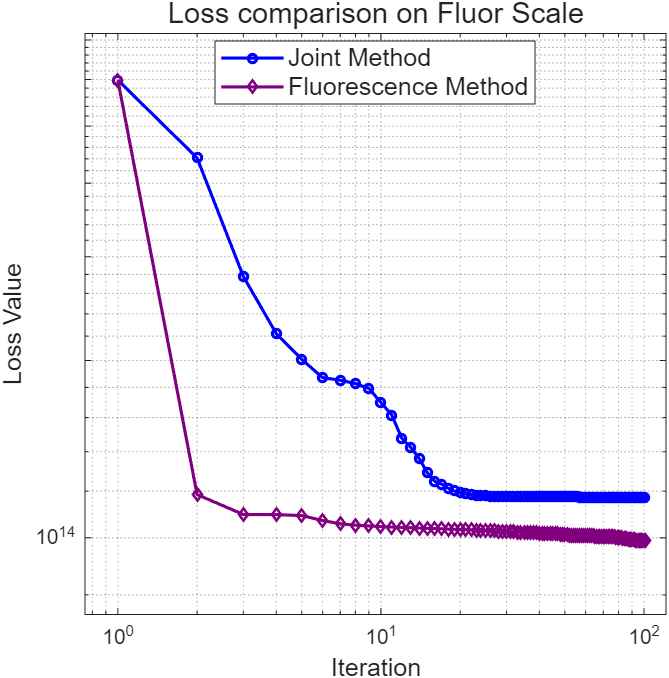}
    \vspace{0.6em}
    \includegraphics[width=0.8\linewidth]{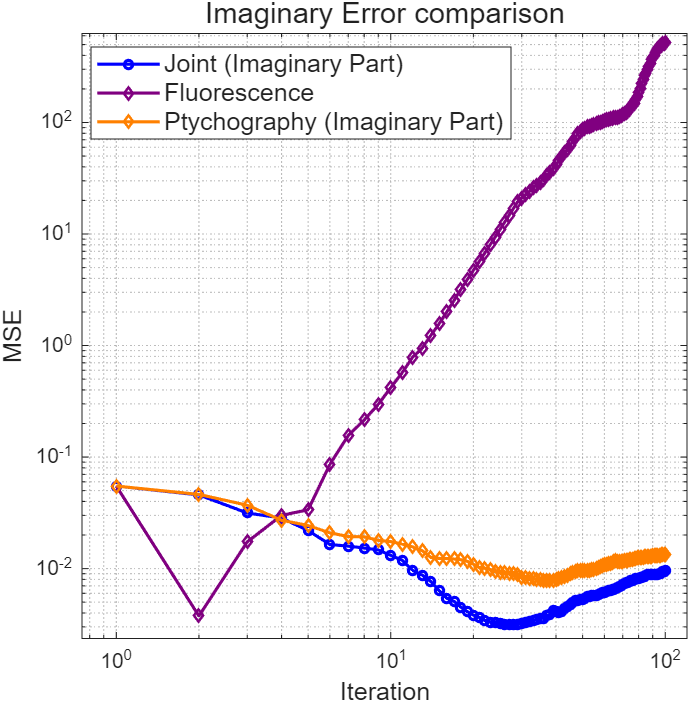}
    \caption*{\(\text{Noise Level} = 10\%\)}
  \end{subfigure}

  \caption{Log--log plots of loss (top) and reconstruction error (bottom) for fluorescence and joint method applied to the synthetic object shown in Figure \ref{fig:noise_recon} with \(n=334\), probe size \(m=64\), and 53\% overlap. The joint method has a relatively lower MSE in the high-noise setting.}
  \label{fig:loss_fluor_errro_noise_compare}
\end{figure}

Figure~\ref{fig:noise_recon} shows reconstruction results for both settings, while Figures~\ref{fig:loss_errro_noise_compare}, \ref{fig:loss_fluor_errro_noise_compare} present the corresponding loss and error curves. Across both noise levels, the joint method consistently outperforms ptychography-only method. In the noise-free case, this robustness arises because the fluorescence modality provides a strong regularization effect on the imaginary part of the reconstruction, thereby preserving resolution even when the ptychography measurements are corrupted on the edges of the image. In the high-noise setting, the advantage in noise robustness becomes more pronounced. The strong regularization provided by ptychography prevents the joint reconstruction from being corrupted by the input noise, unlike the {fluorescence-only} reconstruction.

\subsection{Out-of-focus and in-focus probe}

{
We consider probes with different focal distances, corresponding to in-focus and out-of-focus illumination conditions. Under the in-focus probe setting, the specimen is positioned near the focal plane, resulting in a relatively small and tightly confined illumination spot on the sample. This spatially localized in-focus probe can be regarded as a more tightly confined illumination function, which may reduce the effective overlap ratio in ptychography and thereby further exacerbate the ill-posedness of the ptychographic reconstruction problem. On the other hand, a broader or more defocused probe produces stronger fluorescence blurring and a more ill-conditioned deconvolution problem, while a narrower or more focused probe yields higher intrinsic XRF resolution. As shown in Figure \ref{fig:ptycho_probe_setting}, the mean squared error gap between the ptychographic-only and joint methods increases in the in-focus setting, reflecting the greater ill-posedness of the reconstruction problem in this regime. In addition, as shown in Figure \ref{fig:fluor_probe_setting}, the more focused probe provides higher intrinsic XRF resolution, allowing the fluorescence-only reconstruction to attain a relatively lower MSE during the initial iterations compared with the out-of-focus probe setting. }

\begin{figure}
  \centering
  \begin{subfigure}[t]{0.45\linewidth}
    \centering
    \includegraphics[width=\linewidth]{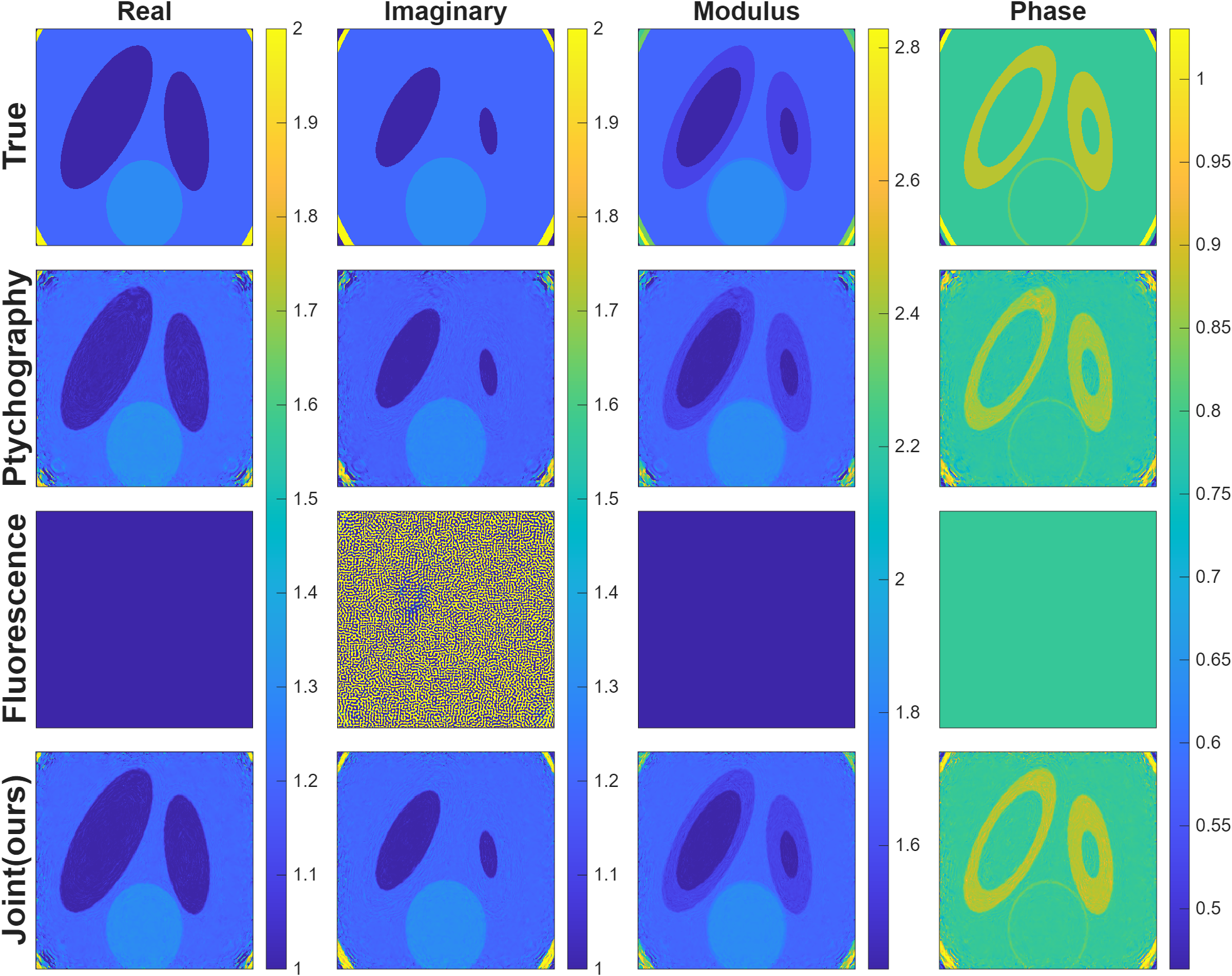}
    \caption*{Out-of-focus probe}
  \end{subfigure}%
  \hfill
  \begin{subfigure}[t]{0.45\linewidth}
    \centering
    \includegraphics[width=\linewidth]{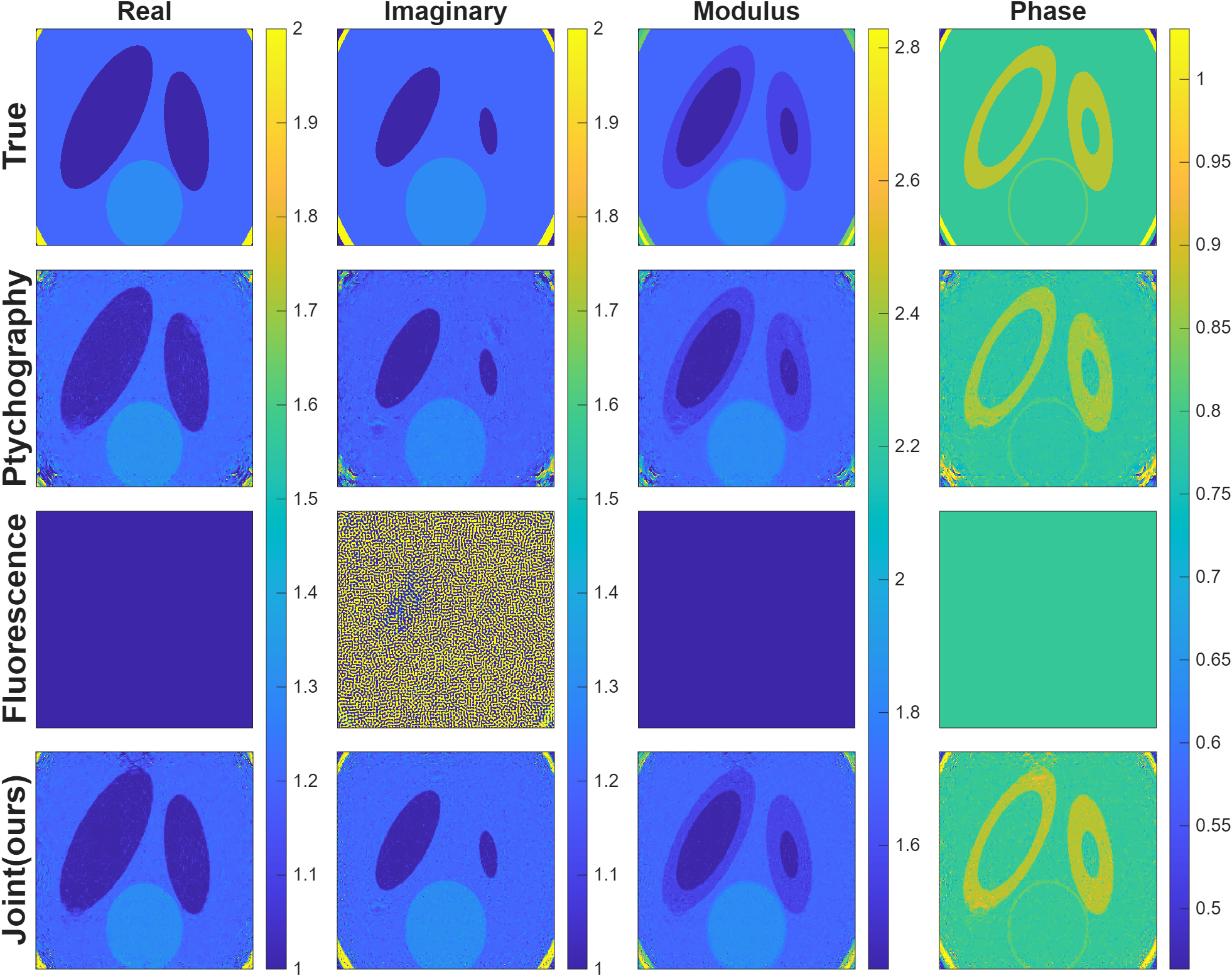}
    \caption*{In-focus probe}
  \end{subfigure}

  \caption{{ Reconstruction results for ptychography, fluorescence, and the joint method applied to a synthetic object under different probe conditions: out-of-focus (left) and in-focus (right).}}
  \label{fig:recon_probe_setting}
\end{figure}

\begin{figure}{ht}
  \centering
  \begin{subfigure}[t]{0.45\linewidth}
    \centering
    \includegraphics[width=0.8\linewidth]{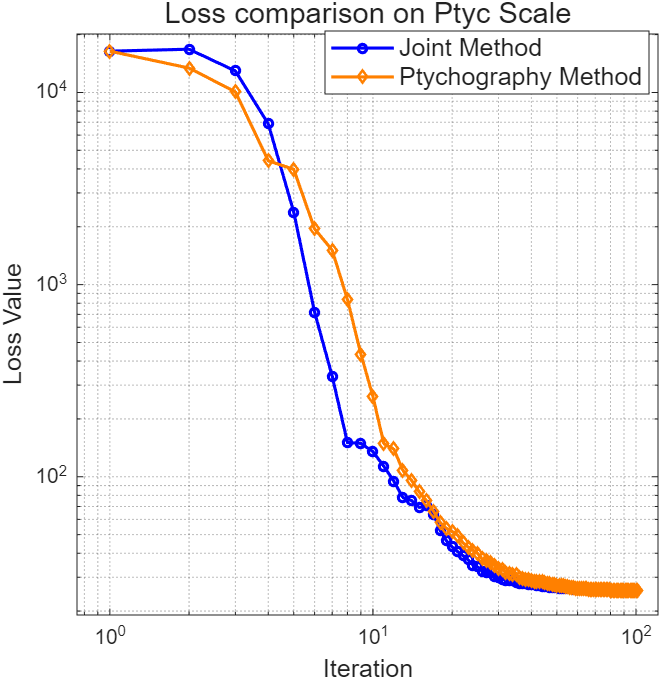}
    \vspace{0.6em}
    \includegraphics[width=0.8\linewidth]{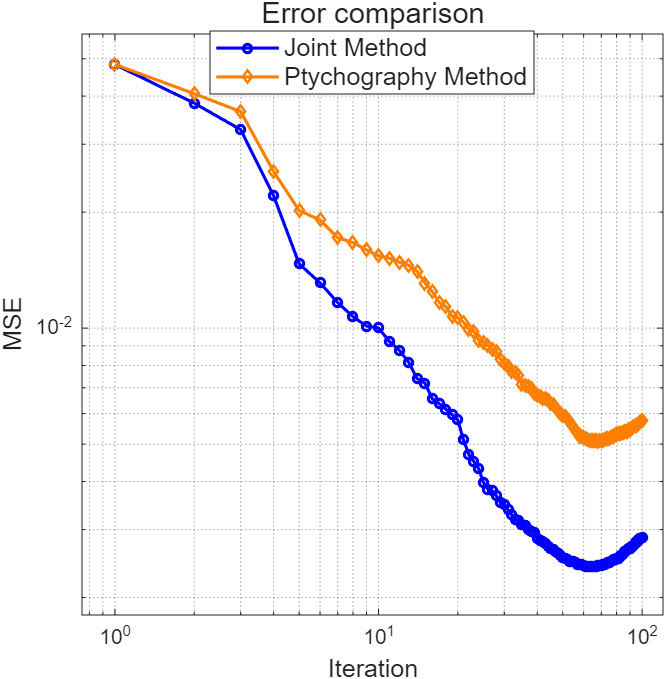}
    \caption*{Out-of-focus probe}
  \end{subfigure}%
  \hfill
  \begin{subfigure}[t]{0.45\linewidth}
    \centering
    \includegraphics[width=0.8\linewidth]{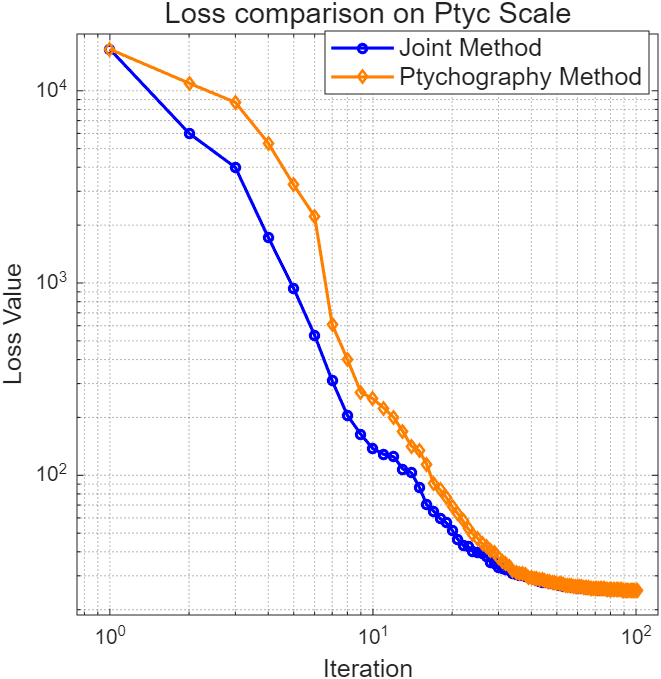}
    \vspace{0.6em}
    \includegraphics[width=0.8\linewidth]{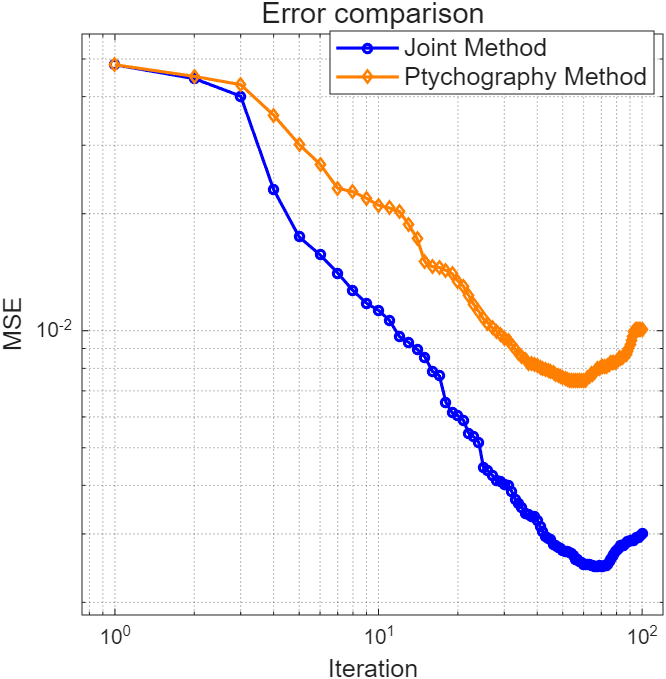}
    \caption*{In-focus probe}
  \end{subfigure}

  \caption{{Log–log plots of the objective loss (top) and reconstruction error (bottom) for ptychography and the joint method under different probe settings for the object shown in Figure \ref{fig:recon_probe_setting}. Under the in-focus condition, the joint method produces a larger MSE gap relative to the ptychography method. }}
  \label{fig:ptycho_probe_setting}
\end{figure}

\begin{figure}{ht}
  \centering
  \begin{subfigure}[t]{0.45\linewidth}
    \centering
    \includegraphics[width=0.8\linewidth]{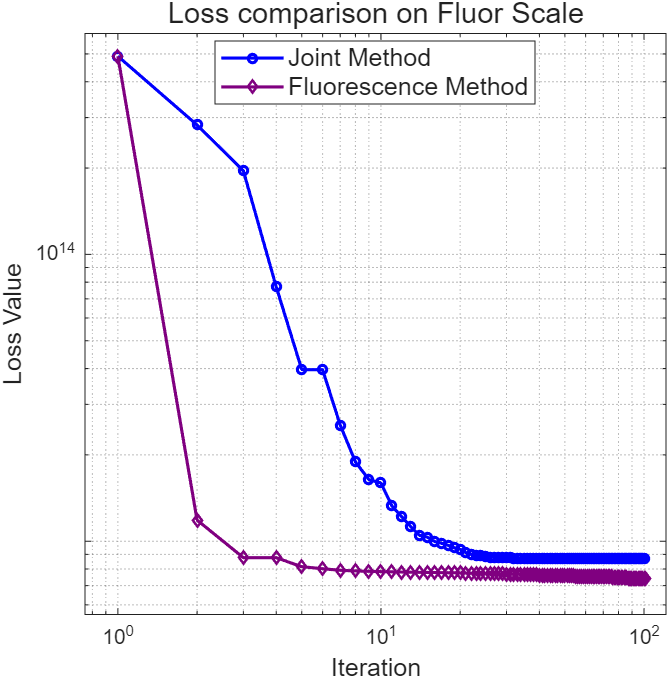}
    \vspace{0.6em}
    \includegraphics[width=0.8\linewidth]{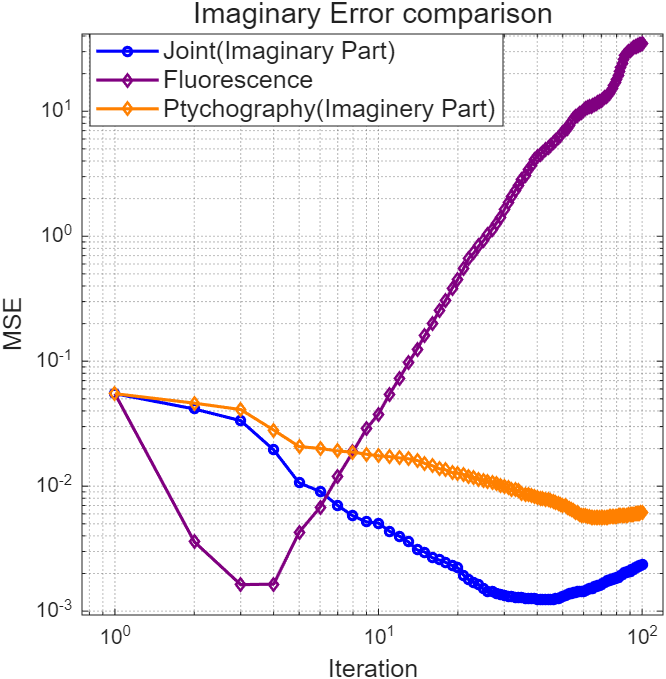}
    \caption*{Out-of-focus probe}
  \end{subfigure}%
  \hfill
  \begin{subfigure}[t]{0.45\linewidth}
    \centering
    \includegraphics[width=0.8\linewidth]{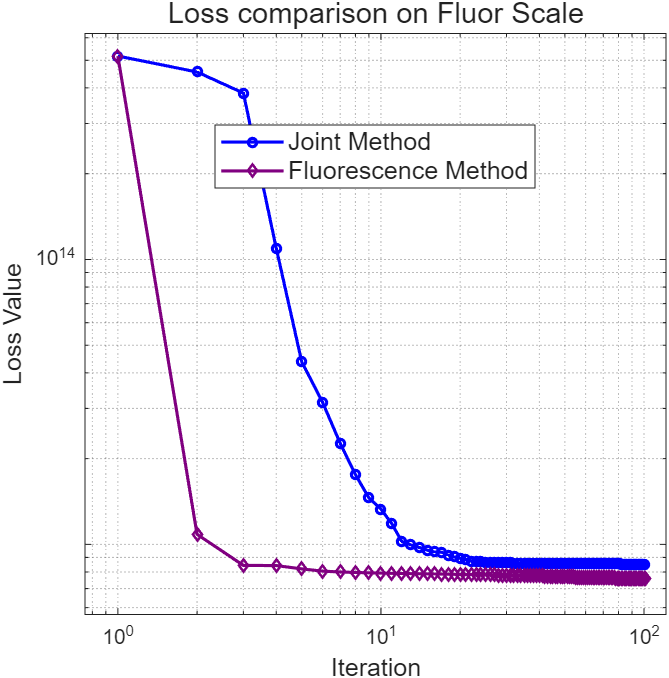}
    \vspace{0.6em}
    \includegraphics[width=0.8\linewidth]{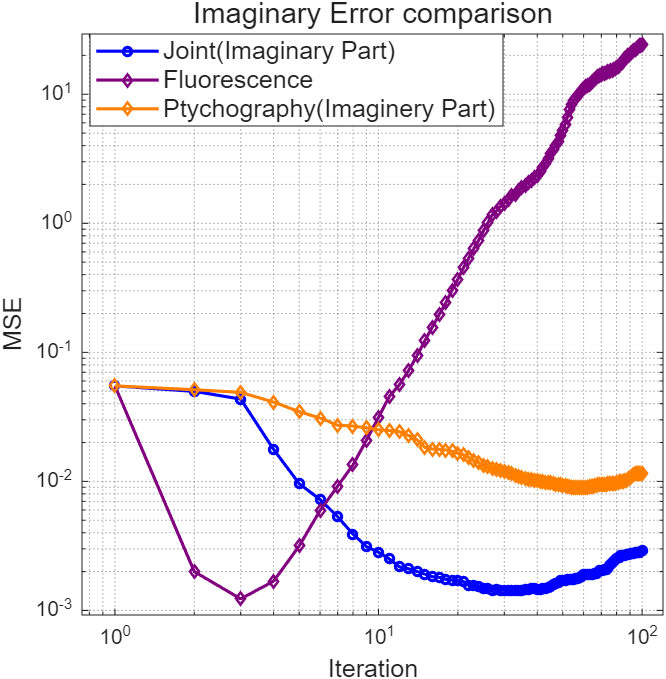}
    \caption*{In-focus probe}
  \end{subfigure}

  \caption{{Log--log plots of objective loss (top) and reconstruction error (bottom) for fluorescence and the joint method under different probe settings for the object shown in Figure \ref{fig:recon_probe_setting}. }  }
  \label{fig:fluor_probe_setting}
\end{figure}

\subsection{Multiple element maps}

To assess the robustness of our joint framework with richer multimodal data, we extend the previous experiment by incorporating multiple elemental concentration maps into the imaginary part of both the joint and ptychography reconstructions. Specifically, we simulate \(N_e = 3\) elements, corresponding to fluorescence measurements \(\{D_e\}_{e=1}^3\).  

\begin{figure}
  \centering
  \begin{subfigure}[t]{0.6\linewidth}
    \centering
    \includegraphics[width=\linewidth]{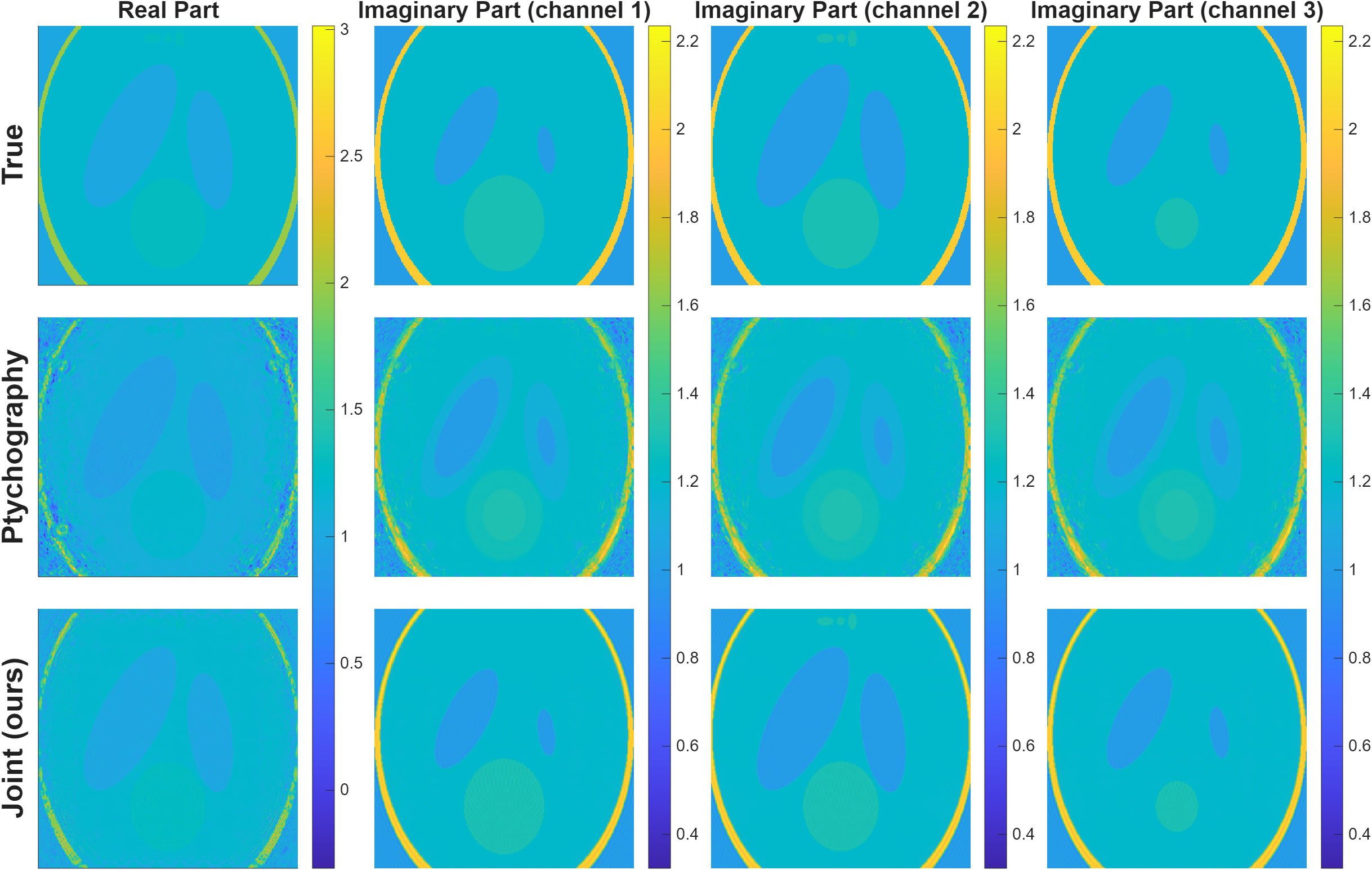}
  \end{subfigure}%
  \hfill
  \begin{subfigure}[t]{0.35\linewidth}
    \centering
    \includegraphics[width=\linewidth]{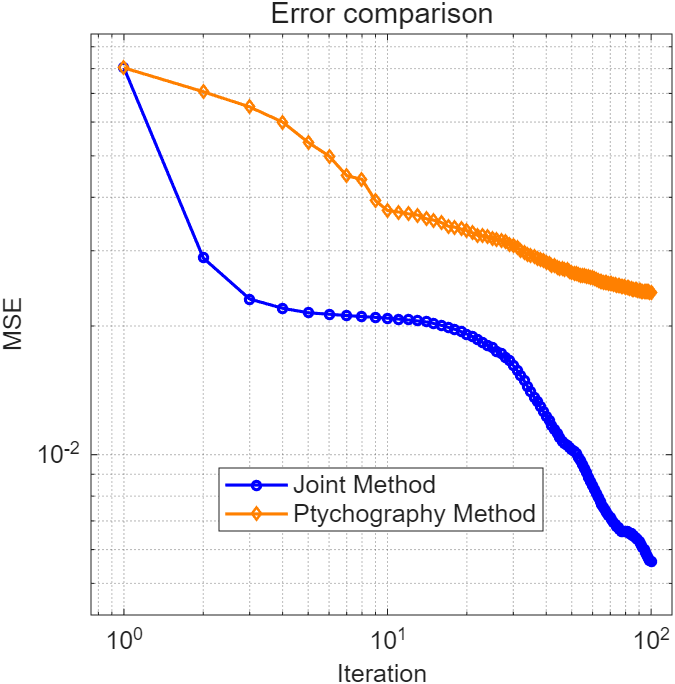}
  \end{subfigure}

  \caption{Joint and ptychography reconstructions running \(200\) iterations on synthetic objects with \(n = 334\), overlap ratio \(= 0.53\), \(\text{Noise Level}= 0\%\), and \(N_e = 3\). Left: reconstruction images across multiple channels. Right: corresponding MSE comparison.}
  \label{fig:recon_img_multiple_ne}
\end{figure}

With the additional information, the joint reconstruction successfully recovers all three elemental maps with better accuracy, as shown in Figure~\ref{fig:recon_img_multiple_ne}. Moreover, the availability of multiple fluorescence channels strengthens the overall constraints on the reconstruction problem: even the \emph{real part} of the joint reconstruction significantly outperforms ptychography, underscoring the benefit of integrating multimodal fluorescence signals as complementary information.  

\subsection{Robustness across different datasets}

\begin{figure}
    \centering
    \begin{subfigure}[t]{0.48\linewidth}
        \centering
        \includegraphics[width=\linewidth]{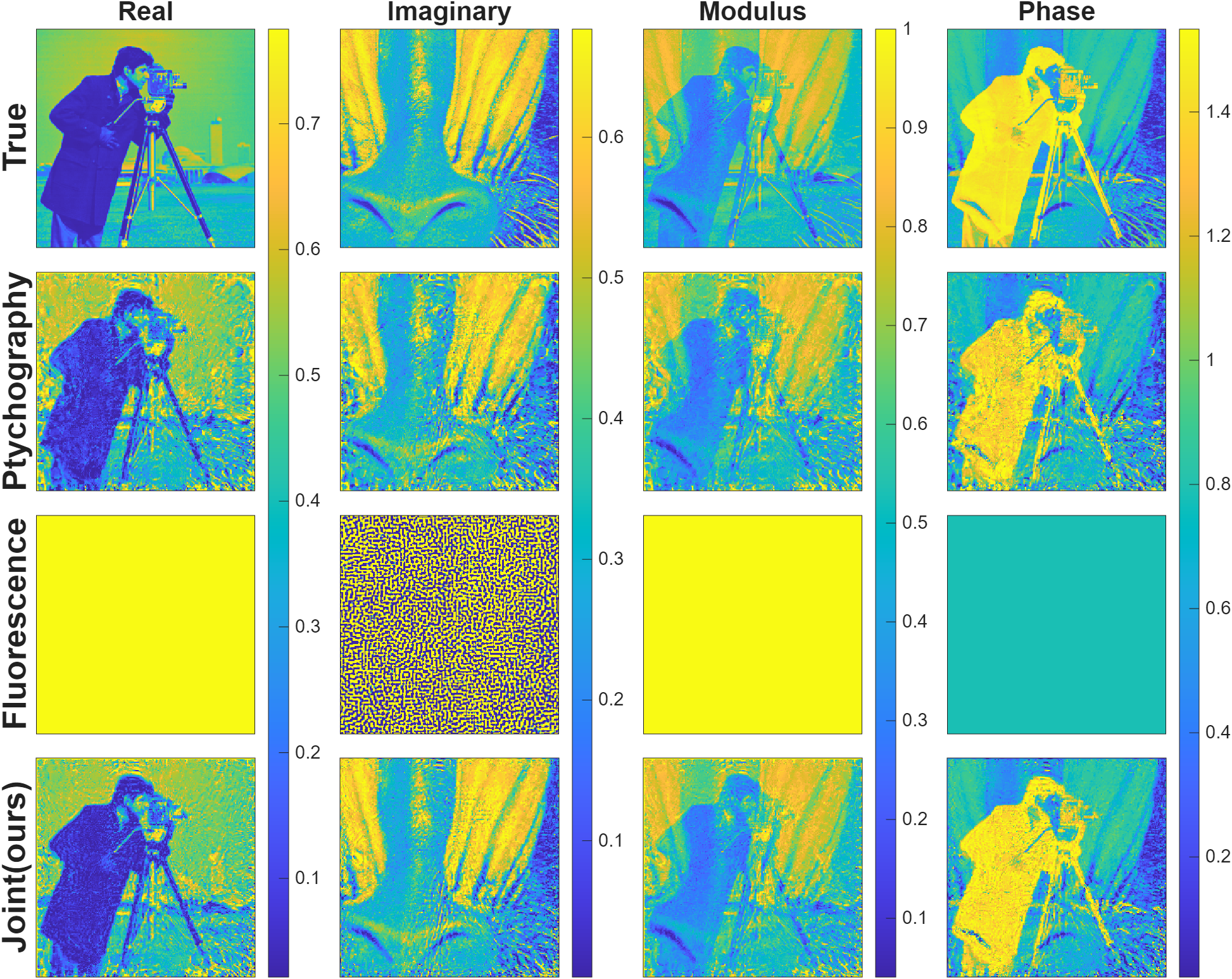}
        \caption*{$n=216$, $m=36$, 44.4\% overlap}
    \end{subfigure}\hfill
    \begin{subfigure}[t]{0.48\linewidth}
        \centering
        \includegraphics[width=\linewidth]{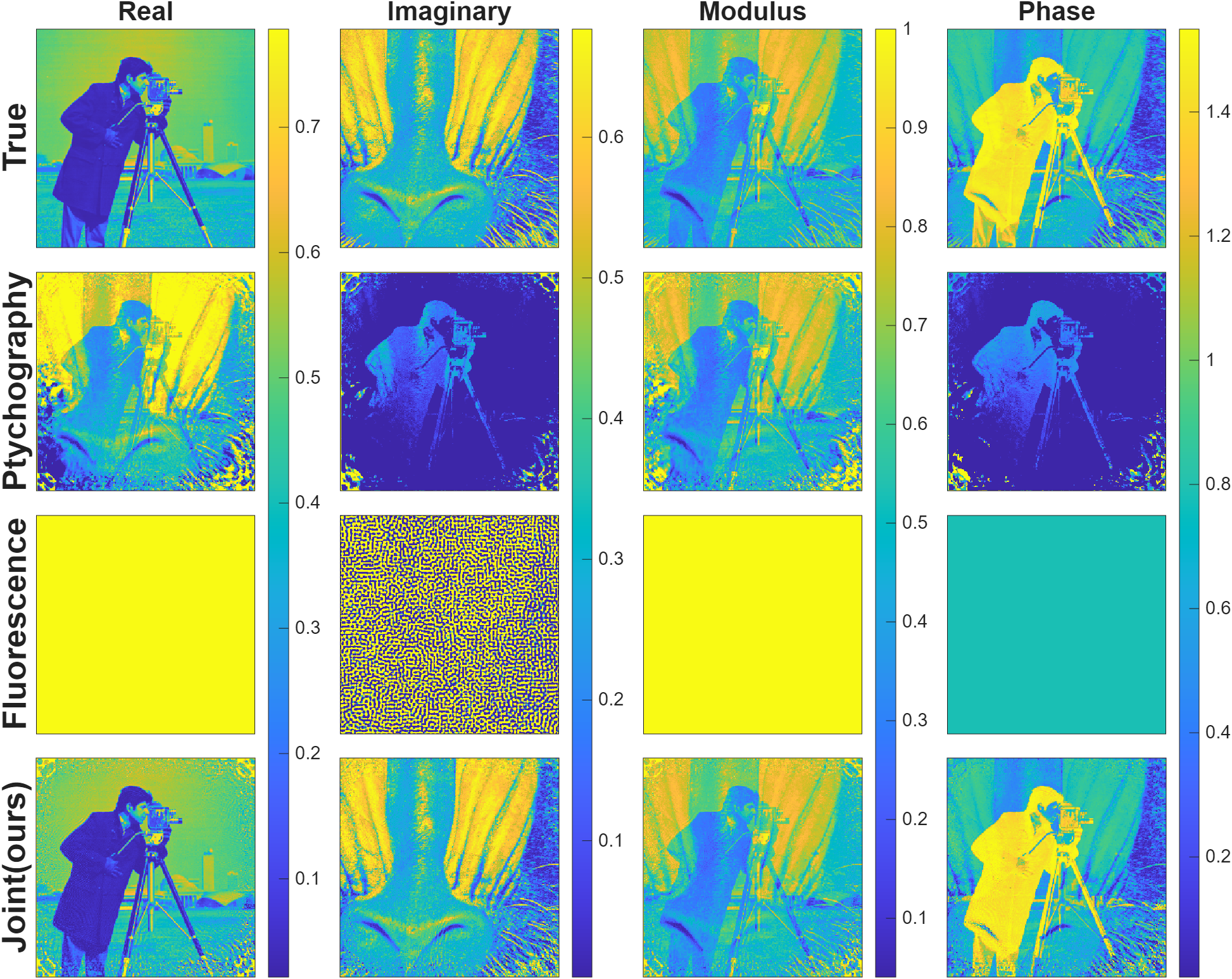}
        \caption*{$n=244$, $m=64$, 68.5\% overlap}
    \end{subfigure}

    \caption{Reconstruction results for the Cameraman--Baboon dataset under $\text{Noise Level} = 3\%$ noise. The joint method produces sharper and more accurate reconstructions than ptychography alone across different parameter settings.}
    \label{fig:recon_cam_baboon_combined}
\end{figure}

\begin{figure}
  \centering
  \begin{subfigure}[t]{0.45\linewidth}
    \centering
    \includegraphics[width=0.8\linewidth]{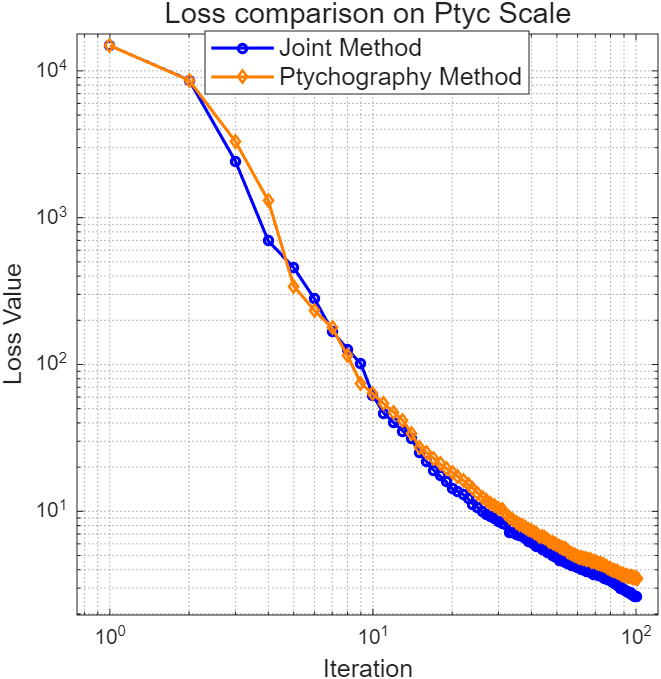}
    \vspace{0.6em}
    \includegraphics[width=0.8\linewidth]{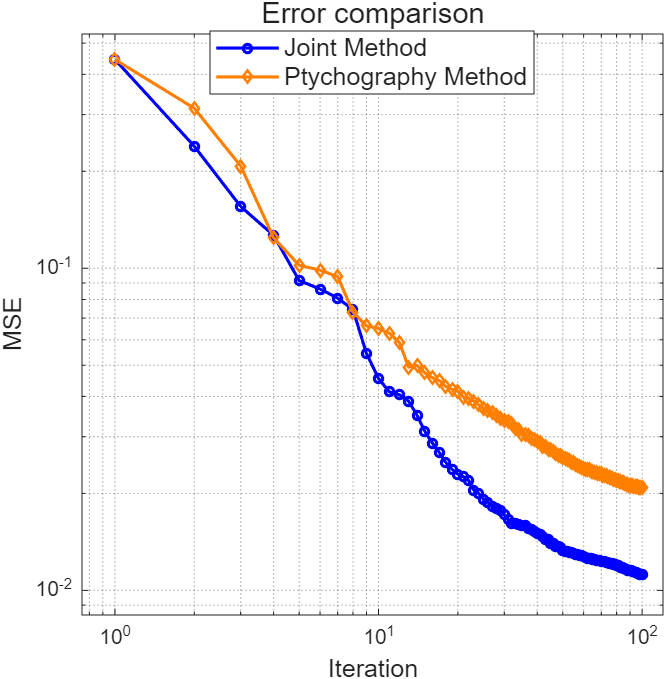}
    \caption*{$n=216$, $m=36$, 44.4\% overlap}
  \end{subfigure}%
  \hfill
  \begin{subfigure}[t]{0.45\linewidth}
    \centering
    \includegraphics[width=0.8\linewidth]{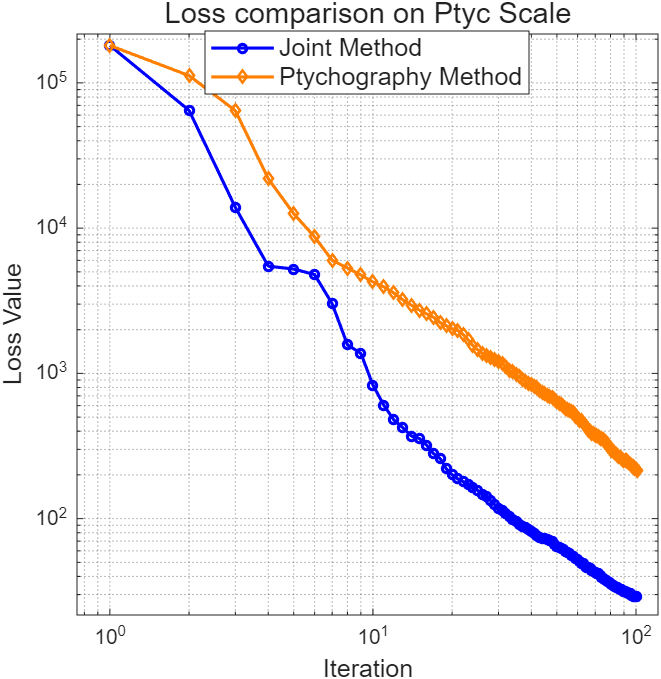}
    \vspace{0.6em}
    \includegraphics[width=0.8\linewidth]{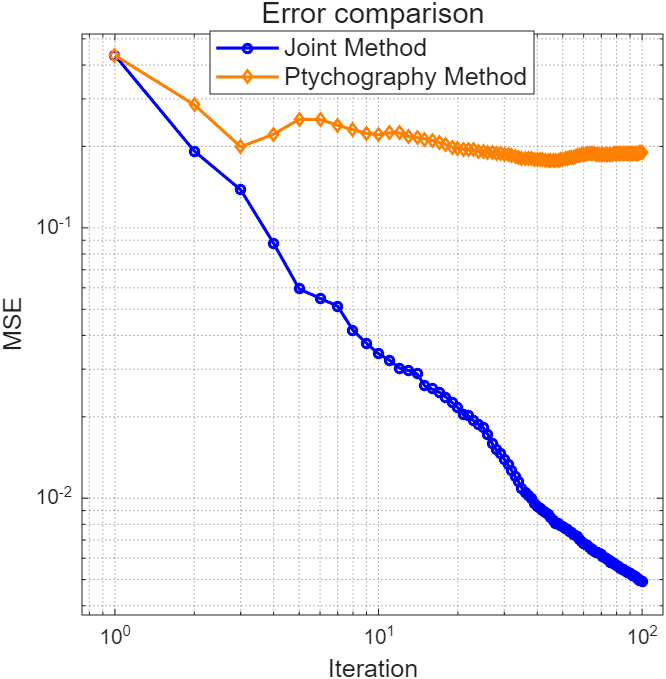}
    \caption*{$n=244$, $m=64$, 68.5\% overlap}
  \end{subfigure}

  \caption{Log--log plots of objective loss (top row) and reconstruction error (bottom row) for ptychography and the joint method on the Cameraman--Baboon dataset shown in Figure \ref{fig:recon_cam_baboon_combined} under $\text{Noise Level} = 3\%$. In both settings, the joint method converges faster and attains lower error.}
  \label{fig:cam_baboon_loss_error_combined}
\end{figure}

\begin{figure}[htbp]
  \centering
  \begin{subfigure}[t]{0.45\linewidth}
    \centering
    \includegraphics[width=0.8\linewidth]{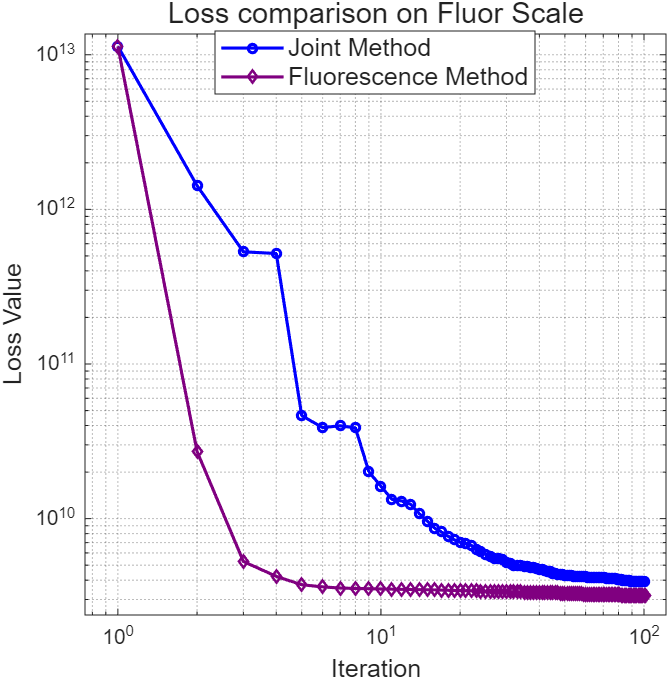}
    \vspace{0.6em}
    \includegraphics[width=0.8\linewidth]{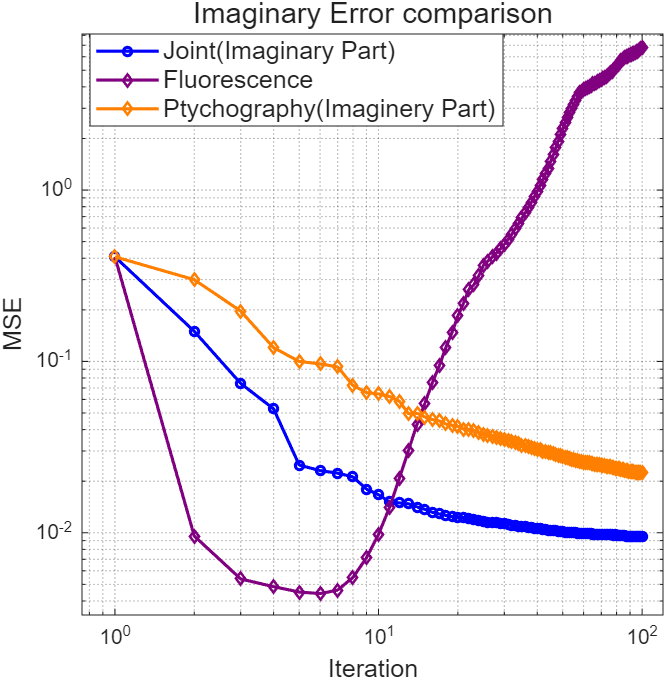}
    \caption*{$n=216$, $m=36$, 44.4\% overlap}
  \end{subfigure}%
  \hfill
  \begin{subfigure}[t]{0.45\linewidth}
    \centering
    \includegraphics[width=0.8\linewidth]{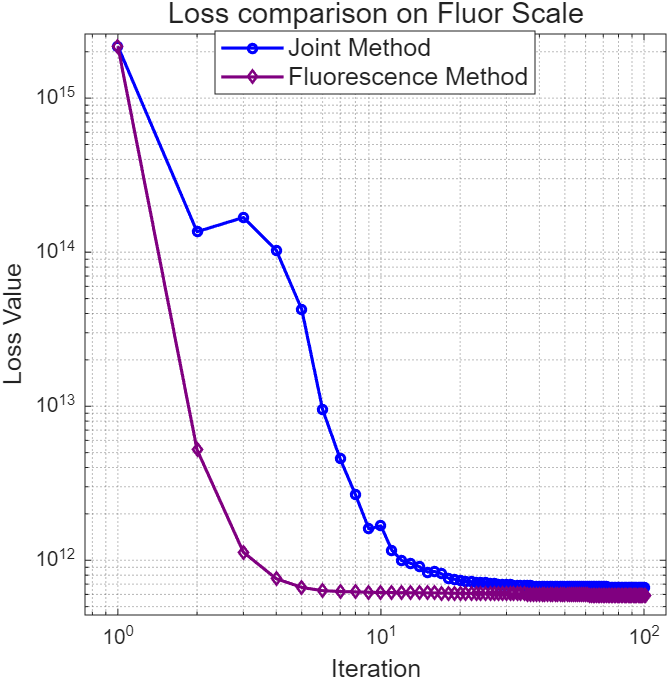}
    \vspace{0.6em}
    \includegraphics[width=0.8\linewidth]{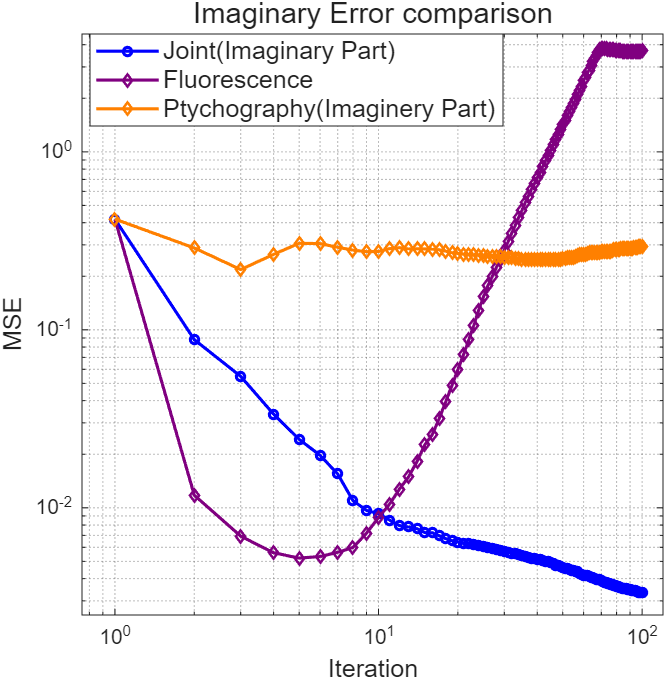}
    \caption*{$n=244$, $m=64$, 68.5\% overlap}
  \end{subfigure}

  \caption{Log--log plots of objective loss (top row) and reconstruction error (bottom row) for fluorescence and the joint method on the Cameraman--Baboon dataset shown in Figure \ref{fig:recon_cam_baboon_combined} under $\text{Noise Level} = 3\%$. }
  \label{fig:fluor_cam_baboon_loss_error_combined}
\end{figure}

We further evaluate the robustness of the joint method on a composite dataset constructed by combining two standard test images: Cameraman (real) and Baboon (imaginary). Both images are normalized to the range $[0,1]$, and the central $256 \times 256$ region of the Baboon image is cropped to match the size of the Cameraman image. Compared with the earlier synthetic dataset, this construction introduces sharper edges and richer textures. As shown in Figures~\ref{fig:recon_cam_baboon_combined} and \ref{fig:cam_baboon_loss_error_combined}, the joint method yields consistently lower errors and clearer reconstructions than ptychography alone. The joint method maintains its superior robustness to input noise compared to fluorescence reconstruction alone, as shown in Figure \ref{fig:fluor_cam_baboon_loss_error_combined}. These results demonstrate that the performance advantage of the joint approach extends across diverse datasets.

\section{Conclusion}

In this work, we present a joint reconstruction framework that integrates ptychography and X-ray fluorescence to enable more robust and quantitative recovery of the object. By combining complementary modalities, the proposed approach mitigates the intrinsic limitations of each method, addressing the resolution constraints of fluorescence imaging and the data incompleteness challenges of ptychography. Numerical experiments demonstrate substantial improvements over single-modality reconstructions, particularly in mean-squared error, highlighting the advantages of exploiting complementary physical information within a unified optimization framework.


The current framework has several important limitations. First, XRF provides chemically specific fluorescence information but does not by itself provide a complete or fully calibrated measurement of absorption. As a result, quantitative recovery may be affected by calibration uncertainty, detector response, self-absorption, spectral overlap, and matrix-dependent attenuation. Second, the present coupling between the ptychographic object and the XRF channels is intentionally simplified and should be viewed as a first-step coupling rather than a complete physical model of the phase–attenuation relationship.

From a computational perspective, the joint reconstruction becomes substantially more demanding for realistic experimental data, since one must optimize over a high-resolution complex object together with fluorescence-related variables under a larger number of scan positions and more complicated forward models. A fair cost comparison should account for the fact that the joint method evaluates both modalities per iteration. Since the task ultimately requires reconstruction of both the ptychographic object and the fluorescence image, the appropriate baseline is the combined cost of running the two single-modality reconstructions separately. Under this accounting, the joint method does not increase the total effort required to recover both quantities. Achieving robustness in this setting will likely require more efficient numerical solvers, improved regularization, and careful implementation for large-scale computation. Future work will focus on testing the method with real experimental data and on extending the model to account for calibration uncertainty, detector effects, and other sources of experimental mismatch.

In addition, the effective XRF resolution is probe-dependent, and the present study does not provide a systematic characterization across probe sizes or focus conditions. Our default experiments use an out-of-focus probe, and we include a limited comparison with an in-focus probe, but a broader study of how probe size and focal condition affect convergence, stability, and reconstruction accuracy remains an important direction for future work. While the present experiments demonstrate the behavior of the joint method under the tested probe settings, a more systematic study of scaling with object size, feature complexity, and probe variation remains an important direction for future work.

\section*{Acknowledgments}
This material was based upon work supported by the U.S.\ Department of Energy, Office of Science,
Office of Advanced Scientific Computing Research's applied mathematics program under contract DE-AC02-06CH11357.

\bibliographystyle{plain}  
\bibliography{references.bib}  

\appendix

\end{document}

%% file: ex_share.tex
\usepackage{graphicx}  
\usepackage{cite}      
\usepackage{subcaption}
\usepackage{placeins} 
\usepackage{amsmath}
\usepackage{amssymb}
\usepackage{array}
\usepackage{xcolor}
\usepackage{float}        
\usepackage{placeins}    
\usepackage{booktabs}
\usepackage{algorithm}
\usepackage{algorithmic}
\usepackage{svg}
\usepackage{placeins}

\headers{Joint X-ray Reconstruction for Multimodal Imaging}{E. Zou, E. Buser, Z. W. Di, Y. Xi}

\newcommand{\WDNote}[1]           
{\textcolor{red}{#1}\marginpar{\textcolor{red}{WD $\longleftarrow$}}}

\author{Chengru Eric Zou\thanks{Department of Mathematics, Emory University, Atlanta, 30322 (\email{eric.zou, ebuser, yxi26@emory.edu})} The research of C. Zou and Y. Xi is supported by NSF awards DMS-2338904 and DMS-2038118.
\and Elle Buser\footnotemark[1]
\and Zichao Wendy Di\thanks{Mathematics and Computer Science Division \& Advanced Photon Science, Argonne National Lab (\email{wendydi@anl.gov})}
\and Yuanzhe Xi\footnotemark[1]
}